\documentclass[12pt]{amsart}

\usepackage[letterpaper, margin=1in]{geometry}
\usepackage{setspace}

\usepackage{verbatim}
\usepackage{amsmath,amsthm,amsfonts,amssymb,mathrsfs}
\numberwithin{equation}{section}
\usepackage{xcolor} 
\usepackage{graphicx}
\usepackage{enumitem}
\usepackage{tikz}
\usepackage{tikz-cd}
\usepackage{pgfplots}
\pgfplotsset{compat=newest}  
\usepackage[colorlinks=true, linkcolor=blue, citecolor=blue, urlcolor=blue]{hyperref}

\DeclareMathOperator{\Ric}{Ric}

\DeclareMathOperator{\im}{Image}

\DeclareMathOperator{\foc}{foc}
\DeclareMathOperator{\cut}{cut}

\DeclareMathOperator{\diam}{diam}
\DeclareMathOperator{\dist}{dist}
\DeclareMathOperator{\tr}{tr}
\DeclareMathOperator{\vol}{vol}
\DeclareMathOperator{\area}{area}

\DeclareMathOperator{\Lip}{Lip}
\DeclareMathOperator{\grad}{grad}

\newcommand{\R}{\mathbb{R}}

\newtheorem{theorem}{Theorem}[section]
\newtheorem*{theorem*}{Theorem}
\newtheorem{proposition}{Proposition}[section]

\newtheorem{corollary}{Corollary}[section]

\newtheorem{lemma}{Lemma}[section]
\newtheorem{definition}{Definition}[section]
\theoremstyle{remark}
\newtheorem{remark}{Remark}[section]
\theoremstyle{definition}
\newtheorem{example}{Example}[section]

\begin{document}
\begin{abstract}
	In this paper, we prove $W^{1,p}$ ($p>n$) and $C^{0,\alpha}$ ($0 < \alpha < 1$) precompactness for classes of Riemannian $n$-manifolds with boundary satisfying uniform $L^{\infty}$
	bounds on curvature, mean curvature, diameter, and the $(n-1)$-volume of the boundary. In particular, we identify a class of convex manifolds and a class of uniformly mean-convex manifolds that
	exhibit this compactness without any a priori bounds on the injectivity radii. In these cases, we show the boundary data excludes interior volume collapse.

	If injectivity radius bounds are assumed, the sectional curvature bounds can be replaced with bounds on the Ricci tensor. A central feature of these results is that we only assume a
	pointwise bound on the mean curvature, which constrains the regularity of the convergence.

	We obtain two geometric stability theorems as applications of this compactness result, based
	on rigidity results of Cohn-Vossen/Pogorelov and Hopf. The first theorem applies to $3$-manifolds with pointwise curvature close to $0$ whose boundaries are Gromov-Hausdorff close to a
	fixed metric on $S^2$ with positive curvature. We show that
	such manifolds are $C^{0, \alpha}$ close to the solid region enclosed by a Weyl embedding of the fixed metric into $\R^3$.

	The second theorem demonstrates that if a $3$-manifold with
	$\chi(\partial M) = 2$ has Ricci curvature close to $0$ and boundary mean curvature close to $1$, then $M$ is $C^{0, \alpha}$ close to the Euclidean unit ball in $\R^3$. In these theorems, the geometric stability is within a precompact class identified above, so that the  notion of ``close'' in the hypotheses depends on the constants determining that class.
\end{abstract}

\title[Compactness Theorem and Applications]{A
	Compactness Theorem for Riemannian Manifolds with Boundary and Applications}
\author[]{Kenneth S. Knox}
\address{Department of Mathematics, Stony Brook University, Stony Brook, NY 11794}
\address{Division of Science, Math, and Computing, Simon's Rock at Bard College, Barrytown, NY 12507}
\subjclass[2020]{Primary 53C21, 58J32; Secondary 35J57, 53C45}

\keywords{Riemannian manifolds with boundary, geometric compactness theorems, boundary harmonic radius, mean-convex boundaries, weak Neumann problem, geometric stability, Pogorelov rigidity.}

\maketitle

\tableofcontents

\section[Introduction]{Introduction}

The purpose of this paper is to establish compactness theorems for Riemannian manifolds with boundary under $L^{\infty}$ curvature bounds and to demonstrate some geometric stability theorems as examples of this theory. The only example of previous work we know, assuming only $L^{\infty}$ curvature bounds, is due to Kodani \cite{MR1065204} and is based solely on comparison geometry applied to a class of convex manifolds. Specifically, for a fixed dimension $n$ and positive constants $K_0, S_0, D_0$, and $v_0$, write $\mathcal{M}_c = \mathcal{M}_c(n, K_0, S_0, D_0, v_0)$ for the class of smooth, compact, connected Riemannian $n$-manifolds with boundary satisfying
\begin{gather*}
	|\sec(M)| \leq K_0, \quad 0 \leq S \leq S_0, \\
	\diam(M) \leq D_0, \quad \vol(M) \geq v_0.
\end{gather*}
Here the bound on $\sec(M)$ means that for each $p \in M$ and each 2-plane $\Pi \in T_pM$, the sectional curvature satisfies $|\sec(\Pi)| \leq K_0$. The quantities $\diam(M)$ and $\vol(M)$ refer to the diameter and volume, respectively, and $S$ denotes the second fundamental form of the boundary. Then \cite[Theorem B]{MR1065204} asserts that any sequence in $\mathcal M_c$ admits a subsequence converging in the Lipschitz topology to a Riemannian manifold with a $C^0$ metric tensor. That is, eventually each term of the sequence is bi-Lipschitz homeomorphic to the limit space with dilatation close to one (cf. also \cite[page 178]{MR2379775} for some further remarks on this work). However, one expects that convergence in a stronger sense is possible. To see why, we note that Anderson showed in \cite{BoundaryMetric3mfld} that
\begin{gather*}
	\Ric(g), H, [h]
\end{gather*}
form an elliptic boundary problem for $g$, where $h$ is the intrinsic boundary metric, $[h]$ represents the pointwise conformal class of $h$, $\Ric(g)$ is the Ricci tensor, and $H$ is the mean curvature of the boundary. These data are certainly controlled by the constants defining $\mathcal M_c$. Further, in \cite{MR2096795}, Anderson, Katsuda, Kurylev, Lassas, and Taylor consider a class of Riemannian manifolds with boundary, which we will denote as $\mathcal A := \mathcal{A}(n, K_0, i_0, D_0, H_0)$, satisfying
\begin{gather*}
	\|\Ric_M\|_{L^{\infty}(M)} \leq K_0, \quad \|\Ric_{\partial M}\|_{L^{\infty}(\partial M)} \leq K_0, \\
	\diam(M) \leq D_0, \quad \|H\|_{\Lip} \leq H_0, \\
	i_M \geq i_0, \quad i_{\partial M} \geq i_0, \quad i_b \geq i_0,
\end{gather*}
and show that $\mathcal A$ is precompact in the $C^{1,\alpha}$ topology, any $0 < \alpha < 1$. Moreover, to any sequence $M_i \in \mathcal A$, the subsequential limit is in the Zygmund space $C^2_*$. Here, the $C^{1,\alpha}$ topology refers to the metrics themselves converging as $C^{1,\alpha}$ functions, where $C^{1,\alpha}(M)$ refers to the standard H\"older space. The quantities $i_M$, $i_{\partial M}$, and $i_b$ refer respectively to the injectivity radius of $M$, $\partial M$ and the boundary injectivity radius. If the conditions on $\mathcal M_c$ are supplemented with $\|H\|_{\Lip} \leq H_0$, call this class $\mathcal M_c^{\Lip}$, we note that (see Proposition \ref{thm:kodani_improvement}) $\mathcal M_c^{\Lip} \subset \mathcal A$ and can therefore conclude this stronger precompactness result on $\mathcal M^{\Lip}_c$. However, it is not satisfying to add this extra condition. Our main result, Theorem \ref{precompact}, is a compactness theorem in the more general setting when $H$ is bounded in $L^{\infty}$. Therefore we define $\mathcal A^{bdd} = \mathcal{A}^{\mathrm{bdd}}(n, K_0, i_0, D_0, H_0)$, to be the class of manifolds satisfying
\begin{gather*}
	\|\Ric_M\|_{L^{\infty}(M)} \leq K_0, \quad \|\Ric_{\partial M}\|_{L^{\infty}(\partial M)} \leq K_0, \\
	\diam(M) \leq D_0, \quad |H| \le H_0, \\
	i_M \geq i_0, \quad i_{\partial M} \geq i_0, \quad i_b \geq i_0.
\end{gather*}
In Theorem \ref{precompact}, we show that $\mathcal A^{bdd}$ is precompact in the $W^{1,p}$, $p > n$, and $C^{0, \alpha}$, $0 < \alpha < 1$, topologies, where $W^{1,p}(M)$ refers to the Sobolev space of one weak derivative in $L^p(M)$. We show in Proposition \ref{thm:kodani_improvement} that $\mathcal M_c \subset \mathcal A^{bdd}$, improving Kodani's result to $W^{1,p}$ metric convergence.

We establish other geometric conditions that imply inclusion in $\mathcal A^{bdd}$. First, define $\mathcal M = \mathcal M(n, K_0, H_0, D_0, A_0, i_0)$ to be the class of compact, connected, Riemannian $n$-manifolds with boundary satisfying
\begin{gather*}
	|\sec(M)| \leq K_0, \quad |\sec(\partial M)| \leq K_0, \quad |H| \leq H_0,\\
	\diam(M) \leq D_0, \quad \vol_{n-1}(\partial M) \geq A_0, \quad i_b \geq i_0.
\end{gather*}

In this case $\vol_{n-1}(\partial M)$ refers to the $(n-1)$-dimensional volume of the boundary. We show in Corollary \ref{cor:strata_inclusion} that $\mathcal M \subset \mathcal A^{bdd}$, thus $\mathcal M$ is precompact in the $W^{1,p}$ and $C^{0,\alpha}$ topologies. This requires establishing a uniform lower bound on $i_M$ for each element of the class. To do this, we need to show volume cannot collapse in the interior (cf. Proposition \ref{prop:interior_volume}). The idea is to use a uniform collar neighborhood of $\partial M$ and
comparison geometry to assert the existence of a cylinder of a definite size, sitting above any boundary point, and with a definite volume bounded below (cf. Lemma \ref{lem:isoperimetric}). Then we localize a relative volume comparison argument with respect to a corresponding cylinder.

We identify a subclass of $\mathcal M$ that does not require an a priori lower bound on the boundary injectivity radius. We define $\mathcal M_+ = \mathcal M_+(n, K_0, H_0, D_0, A_0)$ to be the class of compact, connected, Riemannian $n$-manifolds with boundary satisfying
\begin{gather*}
	|\sec(M)| \leq K_0, \quad |\sec(\partial M)| \leq K_0, \quad 1/H_0 <|H| \leq H_0,\\
	\diam(M) \leq D_0, \quad \vol_{n-1}(\partial M) \geq A_0,
\end{gather*}
where $H_0 > 1$. This is a class of \emph{uniformly mean convex} manifolds with boundary. We note that the uniform positivity of $H$ implies a lower bound on the boundary injectivity radius (cf. Lemma \ref{lem:cutlocuslemma}).

We analogously define $\mathcal{M}^{\Lip}$ and $\mathcal M^{\Lip}_+$ to assume an additional $\|H\|_{\Lip} \leq H_0$ bound on the mean curvature. Of course, this subsumes the pointwise bound in the former case and complements the positive lower bound in the latter. Correspondingly, $\mathcal M_+^{\Lip} \subset \mathcal M^{\Lip} \subset \mathcal A$, so that these classes
have $C^{1, \alpha}$ compactness for any $0<\alpha <1$ and weak $C^2_*$ convergence. The inclusions are summarized in Figure~\ref{fig:class_inclusions}. As an immediate corollary, all of these classes have diffeomorphism finiteness.

\begin{figure}[htpb]
	\centering
	\begin{tikzcd}[row sep=1.5cm, column sep=1cm]
		{\text{\begin{tabular}{c}
					$W^{1,p}$ and $C^{0,\alpha}$ compactness \\
					\small{(Pointwise bounds on $H$)}
				\end{tabular}}}
		&
		\mathcal{M}_+ \arrow[r, hook] &
		\mathcal{M} \arrow[r, hook] &
		\mathcal{A}^{bdd} &
		\mathcal{M}_c \arrow[l, hook']
		\\
		{\text{\begin{tabular}{c}
					$C^{1,\alpha}$ and weak $C^2_*$ compactness \\
					\small{(Lipschitz bounds on $H$)}
				\end{tabular}}}
		&
		\mathcal{M}_+^{\mathrm{Lip}} \arrow[u, hook] \arrow[r, hook] &
		\mathcal{M}^{\mathrm{Lip}} \arrow[u, hook] \arrow[r, hook] &
		\mathcal{A} \arrow[u, hook] &
		\mathcal{M}_c^{\mathrm{Lip}} \arrow[u, hook] \arrow[l, hook']
	\end{tikzcd}
	\caption{Inclusion hierarchy of the classes. Hooked arrows ($\hookrightarrow$) denote subset inclusions pointing from subclasses to their corresponding supersets. The vertical tiers distinguish the geometric limits constrained by pointwise versus Lipschitz bounds on the mean curvature.}
	\label{fig:class_inclusions}
\end{figure}

A central difficulty in obtaining Theorem \ref{precompact} is the weaker $W^{1,p}$ a priori metric regularity, coming from the fact that the mean curvature only involves one derivative of
the metric. In particular, Propositions \ref{lem:local_strong_convergence} and \ref{weaktostrong} show that a sequence of Riemannian manifolds in $\mathcal A^{bdd}$ converging weakly in $W^{1,p}$, any fixed $p > n$, must in fact converge strongly in $W^{1,q}$, any $q > n$. To do this, we first establish $H^1$ convergence, making use of energy estimates for the weak Neumann problem (cf.\ Lemma \ref{prop:neumann_w12_estimate} and Corollary \ref{cor:local_neumann_w12}) together with compact embeddings of the boundary data coming from the dual Sobolev embeddings. We combine this with Sobolev interpolation and $W^{1,q}$ boundary regularity that comes from Mitrea-Taylor layer potential estimates for the weak Neumann problem (cf.\ Proposition \ref{prop:neumann_l_s_estimate}, Corollary \ref{cor:neumann_gradient_trace}, and Corollary \ref{prop:local_neumann_l_s}).

We apply the convergence results to prove two geometric stability theorems regarding $3$-manifolds with boundary. The first can be viewed as a generalization of Cohn-Vossen's~\cite{Cohn-Vossen} and Pogorelov's~\cite{MR0346714} rigidity theorems. Cohn-Vossen's rigidity theorem states the following.

\begin{quote}
	Let $S_1, S_2 \subset \mathbb{R}^3$ be compact, connected, $C^2$ surfaces without boundary, possessing strictly positive Gaussian curvature $K > 0$. If $\phi: S_1 \to S_2$ is a Riemannian isometry with respect to their induced first fundamental forms (meaning $\phi^* I_2 = I_1$), then $\phi$ is the restriction of a rigid Euclidean motion. That is, there exists an isometry $\Phi \in E(3)$ such that $\Phi|_{S_1} = \phi$.
\end{quote}
Pogorelov generalized this as follows.

\begin{quote}\label{thm:pogorelov_intro}
	A \emph{closed convex surface} $S \subset \mathbb{R}^3$ is defined as the boundary of a compact, convex set with non-empty interior. It is equipped with an intrinsic metric $d_S: S \times S \to \mathbb{R}_{\ge 0}$, where $d_S(x, y)$ is the infimum of the lengths of curves on $S$ connecting $x$ and $y$.
	Let $S_1, S_2 \subset \mathbb{R}^3$ be two closed convex surfaces. If there exists an intrinsic isometry $\phi: (S_1, d_{S_1}) \to (S_2, d_{S_2})$, then $S_1$ and $S_2$ are congruent. That is, $\phi$ is the restriction of a rigid Euclidean motion $\Phi \in E(3)$.
\end{quote}

We show that this rigidity is geometrically stable in the class $\mathcal M_+$. Thus, consider a 3-manifold $(M,g) \in \mathcal M_+$ with boundary $(\partial M, h)$. If $(\partial M, h)$ is Gromov-Hausdorff close to a Riemannian surface $(\Sigma, h_{\Sigma})$ with positive curvature, and $\Ric(g) \sim 0$, then $(M,g)$ is close in the $C^{0,\alpha}\cap W^{1,p}$ topology to
a convex set $N$ in Euclidean space, where $N \subset \mathbb{R}^3$ is a region bounded by a Weyl embedding of $\Sigma$ into $\mathbb{R}^3$. See Theorem \ref{T1} and Corollary \ref{T1_eps_delta} for the formal statements and proofs. We emphasize that the smallness conditions on $\Ric(g)$ and the Gromov-Hausdorff distance between $\Sigma$ and $\partial M$ depend on the constants that define the class $\mathcal M_+$. Similar versions can be stated with respect to the classes $\mathcal M$, $\mathcal M_c$ and $\mathcal A^{bdd}$.

Hopf's rigidity theorem \cite{MR0040042} states the following.
\begin{quote}
	Let $\Sigma$ be a closed, orientable, smooth surface, and let $f: \Sigma \to \mathbb{R}^3$ be a $C^3$ immersion. If $\Sigma$ is homeomorphic to the sphere $\mathbb{S}^2$, and the mean curvature $H$ of the immersion is constant, then $f(\Sigma)$ is a standard round sphere.
\end{quote}

We show that this rigidity is geometrically stable in the class $\mathcal M_+$. Suppose $(M,g) \in \mathcal M_+$ is a $3$-manifold with connected boundary $(\partial M, h)$ satisfying $\chi(\partial M) = 2$. Write $H$ for the mean curvature. If $\Ric(g) \sim 0$,  $H \sim 1$, then $(M,g)$ is close in the $C^{0,\alpha}\cap W^{1,p}$ topology to the Euclidean unit ball $B \subset \R^3$. In this case the smallness conditions on $\Ric(g)$ and $|H-1|$ depend on the constants that determine $\mathcal M_+$ as well as the $C^{2,\alpha}$ norm of $h$. As before, similar versions can be stated for other subclasses of $\mathcal A^{bdd}$. See Theorem \ref{conformal_stability} for the formal statement and proof of this result.

We give examples in Section \ref{examples} illustrating various aspects of the theory and applications. We develop preliminary results in Section \ref{defs}. Since subsequential limits of sequences in $\mathcal A^{bdd}$ are $W^{1,p}$ Riemannian manifolds, a significant portion of this section is devoted to analysis directly on these spaces. This includes quantitative estimates for the Dirichlet and Neumann problems for the Laplace-Beltrami operator, including the weak Neumann problem, boundary harmonic coordinates, the introduction of a boundary
harmonic radius and demonstration of its continuity, and standard abstract convergence results. Section~\ref{sec:comparison} develops the comparison geometry needed to show the inclusions $\mathcal M_+ \subset \mathcal M \subset \mathcal A^{bdd}$ and $\mathcal M_c \subset \mathcal A^{bdd}$. The precompactness of $\mathcal A^{bdd}$ is shown in Section~\ref{section_compactness}, and the applications are proved in Section~\ref{proof_of_T1}.

	{\bf Acknowledgement.}
A portion of this work appeared in the author's Ph.D. thesis, completed
at Stony Brook University. The author thanks his advisor Michael Anderson
for his guidance, support, and valuable insights into the field of
geometric analysis.

\section{Examples}\label{examples}
We illustrate some features of the theory with the examples below. First, we show the necessity of the lower $(n-1)$-volume bound in the classes $\mathcal M_+$ and $\mathcal M$.
\begin{example}
	Fix a positive constant $R$. For each integer $k \ge 1$, let $S^1(1/k)$ denote the circle of radius $1/k$, and write $B^2(R)$ for the flat 2-dimensional disk of radius $R$. We define $M_k = S^1(1/k) \times B^2(R)$, with the standard flat product metric
	\begin{equation}
		g_k = dt^2 + dr^2 + r^2 d\theta^2.
	\end{equation}
	Here $t \in [0, 2\pi/k)$ parametrizes the circle $S^1(1/k)$, and $(r, \theta)$ are standard polar coordinates on the disk $B^2(R)$ with $r \in [0, R]$.

	The interior is flat, the boundary is the flat torus $S^1(1/k) \times \partial B^2(R) = S^1(1/k) \times S^1(R)$, and the principal curvatures are $0$ and $1/R$, so the mean curvature is $H_k = \frac{1}{2R}$. The diameter satisfies $\text{diam}(M_k) \le \pi/k + 2R$. The boundary injectivity radius is $i_b(M_k) = R$, while the injectivity radius of the boundary $i(\partial M_k) \to 0$.

	The boundary area and total volume are given by
	\begin{align}
		\text{area}(\partial M_k) & = \left(\frac{2\pi}{k}\right)(2\pi R) = \frac{4\pi^2 R}{k} \to 0,    \\
		\text{vol}(M_k)           & = \left(\frac{2\pi}{k}\right)(\pi R^2) = \frac{2\pi^2 R^2}{k} \to 0.
	\end{align}
	In the Gromov-Hausdorff sense, the sequence of 3-manifolds $M_k$ converges to the 2-dimensional flat disk $B^2(R)$.
\end{example}

Next, we show that some kind of convexity argument or other control on the boundary injectivity radius is necessary.
\begin{example}\label{ex:simple_spheres} Consider $M_k \subset \R^n$, defined as follows. Write $r$ for the radial distance, so $r = \| x\|$, where $x \in \R^n$. For $k > 1$, define
	\[
		M_k = \{x \in \R^n \colon 1 - \frac{1}{k} \le r \le 1 \}.
	\]
	Write $\Sigma_k$ for the boundary component corresponding to $r = 1 - \frac{1}{k}$; the other boundary component, independent of $k$, is the unit sphere, which we denote by $S^n$. We have $H_k = -\frac{k}{k-1}$, where $H_k$ is the mean curvature of $\Sigma_k$, and of course the mean curvature of $S^n$ is $1$. We note that $M_k \to S^2$ in the Gromov-Hausdorff topology and that $\vol(M_k) \to 0$, while the interior curvature is $0$ and the intrinsic boundary curvatures, as well as the second fundamental forms, are uniformly bounded. The diameter of $M_k$ and the injectivity radius of $\partial M_k$ are uniformly bounded below. For any $x \in M_k \setminus \partial M_k,$ we have $\frac{i(x)}{\dist(x,\partial M_k)}$ uniformly bounded below, meaning that there is a uniform lower bound on $i(M_k)$, while the boundary injectivity radius $i_b(M_k) \to 0$. Thus, the uniform lower bound $\frac{1}{H_0} < H$  of $\mathcal M_+$ cannot be replaced with a general bound $|H| \leq H_0$ (and the lower bound on $i_b$ cannot be removed from the definition of $\mathcal M$).
\end{example}

The behavior shown in Example \ref{ex:simple_spheres} can still occur even when the boundary is strictly convex, as the next example shows.
\begin{example}
	Let $M_k = S^2 \times [0, 1/k]$ be equipped with the warped product metric
	\begin{equation*}
		g_k = dt^2 + f_k(t)^2 g_{S^2},
	\end{equation*}
	where $g_{S^2}$ is the standard round metric of radius $1$ on $S^2$, and the warping function is $f_k(t) = 1 - \frac{t}{k} + t^2$.

	The boundary consists of two components: $\Sigma_0 = S^2 \times \{0\}$ and $\Sigma_{1/k} = S^2 \times \{1/k\}$. We compute the mean curvature $H_k$. On $\Sigma_0$, the shape operator is $S(X) = -\nabla_X (\partial_t) = -\frac{f_k'(0)}{f_k(0)} X$. Thus, $H(0) = -2 \frac{-1/k}{1} = \frac{2}{k} > 0$. On $\Sigma_{1/k}$, the shape operator is $S(X) = \frac{f_k'(1/k)}{f_k(1/k)} X$. Thus, $H(1/k) = 2 \frac{1/k}{1} = \frac{2}{k} > 0$.

	Therefore, the boundary is strictly convex. The non-trivial sectional curvatures of $M_k$ are given by the mixed planes and the spherical planes
	\begin{align*}
		K(\partial_t, X) & = -\frac{f_k''(t)}{f_k(t)} = -\frac{2}{1 - t/k + t^2},                            \\
		K(X, Y)          & = \frac{1 - (f_k'(t))^2}{f_k(t)^2} = \frac{1 - (-1/k + 2t)^2}{(1 - t/k + t^2)^2}.
	\end{align*}
	For $t \in [0, 1/k]$, as $k \to \infty$, we have $f_k(t) \to 1$ and $f_k'(t) \to 0$. Consequently, $K(\partial_t, X) \to -2$ and $K(X, Y) \to 1$, meaning the sectional curvatures of the sequence are uniformly bounded.

	The area of the boundary is $\area(\partial M_k) = 4\pi f_k(0)^2 + 4\pi f_k(1/k)^2 = 8\pi$, which is uniformly bounded away from zero. The diameter is bounded by $\pi \max f_k + 1/k \le \pi + 1/k \le 4$.

	The cut distance from the boundary is at most $1/(2k)$, which implies the boundary injectivity radius satisfies $i_b(M_k) \le 1/(2k)$. As $k \to \infty$, $i_b \to 0$ and the interior volume collapses to zero, demonstrating that a uniform lower bound $H \ge 1/H_0 > 0$ is strictly necessary to prevent cut-locus collapse.
\end{example}

Finally, we illustrate why sequences in $\mathcal A^{bdd}$ need not converge in $C^1$ and also provide an example of a sequence satisfying the hypotheses of Theorem \ref{T1}.
\begin{example}
	We construct an explicit sequence of Euclidean domains in $\mathcal{M}_+$, given as perturbations of the unit sphere, which demonstrates that $W^{1,p}(M)$ convergence cannot be improved to $C^{1,\alpha}$ for any $\alpha \ge 0$. Let $(r, \phi, \theta)$ denote standard spherical coordinates in $\mathbb{R}^3$, where $r \ge 0$ is the radial distance, $\phi \in [0, \pi]$ is the polar angle, and $\theta \in [0, 2\pi)$ is the azimuthal angle. Fix a constant $\varepsilon \in (0, 1/2)$. For each integer $k \ge 1$, we define the domain $M_k \subset \mathbb{R}^3$ to be the region bounded by the radial graph $r_k(\phi) = 1 + u_k(\phi)$, where
	\begin{equation}
		u_k(\phi) = \frac{\varepsilon}{k^2} \sin(k \cos\phi).
	\end{equation}
	We equip $M_k$ with the restriction of the standard Euclidean metric $g_{Euc}$ and write $(M_k, g_k)$.
	\begin{figure}[htpb]
		\centering
		\includegraphics[width=0.6\textwidth]{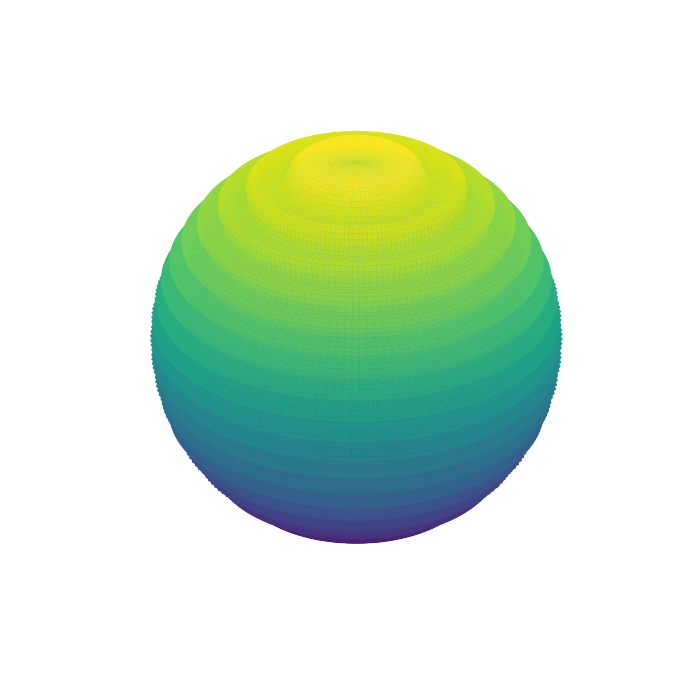}
		\caption{A visualization of $M_{30}$ with $\varepsilon = 36$ to emphasize the perturbation. As $k \to \infty$, $(M_k, g_k)$ converges $W^{1,p}$, $p > n$, to the unit ball with the Euclidean metric, but it does not converge in $C^1$. With $\varepsilon \in (0,\frac{1}{2})$, $M_k$ is eventually convex.}
		\label{fig:wrinkled_sphere}
	\end{figure}
	It is straightforward to verify that $M_k \in \mathcal M_+$ for appropriate constants. We note the $u_k \to 0$ uniformly, and in fact the $W^{1,p}$ geometric limit of $(M_k,g_k)$, guaranteed by Theorem \ref{precompact}, is the unit ball $B$. Let us calculate the mean curvature $H_k$ of $g_k$ as $k \to \infty$. We begin by computing the principal curvatures. Let $\omega(\phi, \theta)$ denote the position vector on $S^2$. The boundary is parametrized by $X(\phi, \theta) = r_k(\phi) \omega(\phi, \theta)$. Therefore the induced boundary metric $h_k$, is diagonal, with components
	\begin{equation}
		(h_k)_{\phi\phi} = r_k^2 + (r_k')^2, \quad (h_k)_{\theta\theta} = r_k^2 \sin^2\phi.
	\end{equation}
	Equipping $\partial M_k$ with the inward-pointing unit normal $N = (r_k^2 + (r_k')^2)^{-1/2} (-r_k \omega + r_k' \omega_\phi)$, it is straightforward to calculate
	\begin{align}
		\kappa_1^{(k)}(\phi) & = 1 + \varepsilon \sin^2\phi \sin(k\cos\phi) + \mathcal{R}_{1,k}(\phi), \\
		\kappa_2^{(k)}(\phi) & = 1 + \mathcal{R}_{2,k}(\phi),
	\end{align}
	where the error bounds satisfy $|\mathcal{R}_{i,k}| \to 0$ as $k \to \infty$. Thus, mean curvature of $\partial M_k$ is
	\begin{equation}
		2H_k(\phi) = 2 + \varepsilon \sin^2\phi \sin(k\cos\phi) + \mathcal{R}_k(\phi),
	\end{equation}
	where $|\mathcal{R}_k| \to 0$. Note that (any subsequence of) $H_k$ deviates from $1$ by at least $\varepsilon/4$ infinitely often as $k \to \infty$. This shows that $
		(M_k, g_k)$ cannot converge in $C^1$ to its limit. In fact, since $\int_{\partial M_k} |1 - H_k|^p \,dS_{h_k}$ is uniformly bounded away from zero for large $k$, this shows that ${g_k}|_{\partial M_k}$ cannot converge strongly in $W^{1,p}(\partial M)$.
	On the other hand, it holds that
	\begin{enumerate}[label=(\roman*)]
		\item $g_k$ converges strongly in $W^{1,q}(B)\cap C^{0,\alpha}(B)$ to $(B, g_{Euc})$, $3 < q < \infty$, $0 < \alpha < 1$,
		\item $h_k$ converges strongly in $C^{1,\alpha}(S^2)$ to the round metric on $S^2$, $0 < \alpha < 1$,
		\item the gradient of $g_k|_{\partial M_k}$ is uniformly bounded in $L^q(\partial M_k)$, any $2 < q < \infty$, and
		\item $H_k$ converges to $1$ weak* on $S^2$.
	\end{enumerate}
	Even though $\partial M_k$ is convex for small $\varepsilon$ and has bounded curvature, control on the extrinsic geometry in $W^{1,p}(\partial M)$ is not a consequence of these constraints.
\end{example}

\section{Preliminary Results} \label{defs}
In this section we set notation, definitions, and preliminary material necessary for our arguments. We begin with some remarks on $W^{1,p}$ Riemannian manifolds, which arise as limits of sequences in $\mathcal A^{bdd}$. We recall some estimates for the Dirichlet problem, show that harmonic coordinate charts exist, and show that the Neumann problem for the Laplacian is well-defined. These are used to demonstrate a boundary regularity result and establish some a priori energy estimates for the Neumann problem. We set definitions for convergence, harmonic and boundary harmonic radii, and show continuity of the boundary harmonic radius. Finally, we note an abstract convergence result which serves to convert a global harmonic radius bounds into weak convergence of manifolds in the corresponding function space.

	{\bf Function Spaces.}
Let us set notation for the function spaces to be used. On a domain $\Omega$, either $\Omega \subset \R^n$ or $\Omega \subset M$, where $M$ is a Riemannian manifold with boundary, write $L^p(\Omega)$ for the space of (equivalence classes of) measurable functions with bounded $L^p$ integral. Write $W^{k,p}(\Omega)$, $k \in \mathbb N$, for the Sobolev space of functions on $\Omega$ with $k$ weak derivatives in $L^p(\Omega)$, and write $W_0^{k,p}(\Omega)$ for the closure of $C^{\infty}_c(\Omega)$ with respect to the $W^{k,p}$ norm. Then $W^{-k,p}(\Omega) := \big (W_0^{k,p'}(\Omega)\big )^*$, the dual space of $W_0^{k,p'}(\Omega)$, where $1/p + 1/p' = 1$. We write $W^{k-1/p, p}(\partial \Omega)$ for the space of traces of $W^{k,p}(\Omega)$ functions, which coincides with the Besov space $B^{k-1/p}_{p,p}(\partial \Omega)$. We write $H^s(\Omega)$ for $W^{s,2}(\Omega)$, typically with $s = 1, -1, 1/2$, or $ -1/2$. There are many excellent references on function space theory, but see for example \cite{MR2759829, MR3024598,MR1163193,MR2250142,MR4298338}.

\subsection{\texorpdfstring{$W^{1,p}$}{W1,p} Riemannian Manifolds}
By a $W^{1,p}$ Riemannian manifold with boundary, we mean a manifold with boundary carrying a $W^{2,p}$ atlas of
charts and a metric tensor $g \in W^{1,p}(M)$. Throughout the discussion we assume $p > n$. Write $h$ for the
induced metric tensor on $\partial M$. By the trace theorem, $h \in W^{1-1/p,p}(\partial M)$. However,  our working hypotheses on $\partial M$ will guarantee $h \in W^{2,p}(\partial M)$ and that $\partial M$ carries an atlas of $W^{2,p}$ charts. The condition $p > n$ implies, via the Sobolev embedding, that $g \in C^{0,\alpha}(M)$, $\alpha = 1 - n/p$. Similarly, the coordinate charts are then $C^{1,\alpha}(M)$ and in particular, are diffeomorphisms. Thus, $M$ carries with it a smooth structure compatible with the $C^1$ atlas (cf. \cite[Theorem 5.11]{MR198479}), and in fact this smooth structure can be taken to be compatible with the $W^{2,p}$ structure on $M$.

As is well-known from \cite{MR644518}, Riemannian manifolds with $C^{k,\alpha}$ metric tensors, $k \geq 1$, admit atlases
of harmonic coordinate charts with transition maps that are at least $C^{k+1,\alpha}$. Let us show a similar result for the class of $W^{1,p}$, $p > n$ Riemannian manifolds with boundary. First, write $\Delta_g$ for the Laplace-Beltrami operator, given in divergence form in coordinates $x^i$ by

\begin{equation}\label{divlap}
	\Delta_g u = \frac{1}{\sqrt{|g|}} \frac{\partial}{\partial x^j} \left( \sqrt{|g|} g^{jk} \frac{\partial u}{\partial x^k} \right),
\end{equation}

and in non-divergence form by

\begin{equation}\label{ndivlap}
	\Delta_g u = g^{jk} \frac{\partial^2 u}{\partial x^j \partial x^k} - g^{jk} \Gamma^m_{jk} \frac{\partial u}{\partial x^m},
\end{equation}
where the Christoffel symbols are
\begin{equation}
	\Gamma^m_{jk} = \frac{1}{2} g^{ml} \left( \frac{\partial g_{lj}}{\partial x^k} + \frac{\partial g_{lk}}{\partial x^j} - \frac{\partial g_{jk}}{\partial x^l} \right).
\end{equation}
From the Sobolev embedding, the coefficients of \eqref{divlap} are at least $C^{0,\alpha}$, $\alpha = 1 - n/p$. Moreover, we see that the Christoffel symbols are $L^p$ functions on their domains, so that the first-order coefficients of \eqref{ndivlap} are in $L^p$.
We recall the definition of Harmonic coordinates.
\begin{definition}[Harmonic Coordinates]
	Suppose $(M,g)$ is a $W^{1,p}$ Riemannian manifold and $U \subset M$ is an open set.
	\begin{enumerate}
		\item If $U \cap \partial M = \emptyset$, a coordinate chart $x:U\to \mathbb R^n$, $x = (x^1, \ldots, x^n)$ is an \emph{interior harmonic coordinate chart} if
		      \[
			      \Delta_g x^i = 0
		      \]
		      for each $1 \leq i \leq n$.

		\item If $U \cap \partial M \neq \emptyset$, then a boundary chart $x: U \to \R^n_+$ is a \emph{boundary harmonic coordinate chart} if each component is harmonic on the interior, \[ \Delta_g x^i = 0 \text{ in } U \backslash \partial M \text{ for } i = 1, \dots, n, \]
		      and if the tangential components are harmonic on the intrinsic boundary, \[ \Delta_h x^{\alpha}|_{\partial M} = 0 \text { for } 1 \leq \alpha \leq n-1. \]
	\end{enumerate}
	In either case, we call $x$ a harmonic coordinate chart.
\end{definition}

We record some a priori estimates and solvability for the Dirichlet problem in our context that are near-immediate consequences \cite[Theorems 9.13 and 9.15]{MR1814364}.

\begin{lemma}[Local Boundary $W^{2,p}$ Estimate] \label{thm:local_boundary_w2p_nonhom}
	Let $\Omega \subset \mathbb{R}^n$ be a domain with a smooth boundary portion $T \subset \partial\Omega$. Suppose $g$ is a metric on $\Omega \cup T$ with components $g_{ij} \in W^{1,p}(\Omega)$ for $p > n$.

	If $u \in W^{2, p}(\Omega)$ is a solution of $\Delta_g u = f$ in $\Omega$ with $f \in L^p(\Omega)$, and $u = \phi$ on $T$ with $\phi \in W^{2-1/p, p}(T)$, then for any domain $V \Subset \Omega \cup T$, we have
	\begin{equation} \label{eq:local_boundary_w2p_nonhom}
		\|u\|_{W^{2,p}(V)} \le C \left( \|f\|_{L^p(\Omega)} + \|\phi\|_{W^{2-1/p, p}(T)} + \|u\|_{L^p(\Omega)} \right),
	\end{equation}
	where $C > 0$ is a constant depending only on $n$, $p$, $T$, $V$, $\Omega$, the ellipticity constants of $g$, and $\|g\|_{W^{1,p}(\Omega)}$.
\end{lemma}

\begin{proof}
	We use the non-divergence form \eqref{ndivlap} for $\Delta_g u$. Because $g \in W^{1,p}(\Omega)$ and $p > n$, the principal coefficients satisfy $g^{ij} \in C^0(\Omega \cup T)$ and the lower-order coefficients satisfy $b^k = -g^{ij} \Gamma^k_{ij} \in L^p(\Omega)$. As explicitly noted at the end of \cite[Section 9.5]{MR1814364}, the proof of the local boundary estimate \cite[Theorem 9.13]{MR1814364} holds good under these weakened conditions for the lower-order coefficients when $p > n$.

	To handle the non-homogeneous boundary data, choose an open set $U$ satisfying $V \Subset U \Subset \Omega \cup T$. Since $\overline{V} \cap T$ is compact, we can cover it by a finite number of boundary coordinate charts $\{B_k\}_{k=1}^N$ such that $B_k \Subset U$ and each $B_k \cap (\Omega \cup T)$ is mapped to an upper half-ball $B^+_{R_k} \subset \mathbb{R}^n_+$, with $B_k \cap T$ corresponding to the flat portion on $\{x^n = 0\}$.

	Choose a smooth partition of unity $\{\eta_k\}_{k=1}^N$ subordinate to $\{B_k\}$ such that $\sum_{k=1}^N \eta_k \equiv 1$ on a relative neighborhood of $\overline{V} \cap T$ in $T$. In each chart, $\eta_k \phi$ is compactly supported on the flat boundary and belongs to $W^{2-1/p, p}(\mathbb{R}^{n-1})$ (extended by zero). Applying \cite[Theorem 2.7.2]{MR3024598}, we obtain a function $w_k \in W^{2,p}(\mathbb{R}^n_+)$ whose trace is $\eta_k \phi$.

	To ensure the extension is compactly supported within the chart, we multiply $w_k$ by an $n$-dimensional interior cutoff $\zeta_k \in C^\infty_c(B^+_{R_k} \cup \{x^n=0\})$ such that $\zeta_k \equiv 1$ on the support of $\eta_k \phi$. The product $\Phi_k = \zeta_k w_k$ belongs to $W^{2,p}(\mathbb{R}^n_+)$, inherits the bound $\|\Phi_k\|_{W^{2,p}(\mathbb{R}^n_+)} \le C \|\eta_k \phi\|_{W^{2-1/p, p}(\mathbb{R}^{n-1})}$, has trace $\eta_k \phi$, and is identically zero near the curved interior boundary of the half-ball.

	Pulling these functions back to $\Omega$ and extending by zero, we define the global extension $\Phi = \sum_{k=1}^N \Phi_k \in W^{2,p}(\Omega)$. By construction, $\Phi = \phi$ on a relative neighborhood of $\overline{V} \cap T$ in $T$, and it satisfies the bound $\|\Phi\|_{W^{2,p}(\Omega)} \le C \|\phi\|_{W^{2-1/p, p}(T)}$.

	Defining $v = u - \Phi$, we have $v \in W^{2,p}(\Omega)$ and $v = 0$ on $\overline{V} \cap T$. The function $v$ satisfies $\Delta_g v = f - \Delta_g \Phi$ in $\Omega$. The product of the $L^p$ lower-order coefficients of $\Delta_g$ and the bounded gradient of $\Phi$ remains in $L^p(\Omega)$, giving $\|\Delta_g \Phi\|_{L^p(\Omega)} \le C \|\Phi\|_{W^{2,p}(\Omega)}$. Applying the local homogeneous boundary estimate \cite[Theorem 9.13]{MR1814364} directly to $v$, we get
	\begin{align*}
		\|v\|_{W^{2,p}(V)} & \le C \left( \|\Delta_g v\|_{L^p(\Omega)} + \|v\|_{L^p(\Omega)} \right)                                                    \\
		                   & \le C \left( \|f\|_{L^p(\Omega)} + \|\Delta_g \Phi\|_{L^p(\Omega)} + \|u\|_{L^p(\Omega)} + \|\Phi\|_{L^p(\Omega)} \right).
	\end{align*}
	Utilizing the triangle inequality and our trace bound on $\Phi$ implies the desired estimate \eqref{eq:local_boundary_w2p_nonhom}.
\end{proof}

\begin{lemma}[Dirichlet problem with $g \in W^{1,p}(M)$] \label{thm:strong_solutions}
	Let $\Omega \subset \mathbb{R}^n$ be a bounded domain with smooth boundary. Suppose $g$ is a metric on $\overline{\Omega}$ with $g_{ij} \in W^{1,p}(\Omega)$ for $p > n$. If $\phi \in W^{2-1/p, p}(\partial \Omega)$, then the Dirichlet problem
	\begin{equation}\label{dirichlet_general}
		\begin{cases}
			\Delta_g u = 0 & \text{ in } \Omega,         \\
			u = \phi       & \text{ on } \partial \Omega
		\end{cases}
	\end{equation}
	admits a unique strong solution $u \in W^{2,p}(\Omega)$. Furthermore, there exists a constant $C > 0$, depending only on $n$, $p$, $\Omega$, the ellipticity constants of $g$, and $\|g\|_{W^{1,p}(\Omega)}$, such that
	\begin{equation}\label{eq:apriori_general}
		\|u\|_{W^{2,p}(\Omega)} \leq C \left(\|\phi\|_{W^{2-1/p, p}(\partial \Omega)} + \|u\|_{L^p(\Omega)}\right).
	\end{equation}
\end{lemma}
\begin{proof}
	We emphasize that this follows immediately from the proof of \cite[Theorem 9.15]{MR1814364}. With $p > n$, the argument holds good with $b^k \in L^p(\Omega)$ (the remark at the end of section \cite[Section 9.5]{MR1814364} applies equally well to Theorem 9.15). Thus, extending $\phi$ to a function $\Phi \in W^{2,p}(\Omega)$ and applying Theorem 9.15 gives the existence of $u$.
\end{proof}

\begin{lemma}[Dirichlet problem with $C^0$ boundary data] \label{thm:strong_solutions_c0}
	Let $\Omega \subset \mathbb{R}^n$ be a bounded domain with smooth boundary. Suppose $g$ is a metric on $\overline{\Omega}$ with $g_{ij} \in W^{1,p}(\Omega)$ for $p > n$. If $\phi \in C^0(\partial \Omega)$, then the Dirichlet problem
	\begin{equation}\label{dirichlet_general_c0}
		\begin{cases}
			\Delta_g u = 0 & \text{ in } \Omega,         \\
			u = \phi       & \text{ on } \partial \Omega
		\end{cases}
	\end{equation}
	admits a unique strong solution $u \in W_{\text{loc}}^{2,p}(\Omega) \cap C^0(\overline{\Omega})$.
\end{lemma}

\begin{proof}
	This \cite[Corollary 9.18]{MR1814364}, with our continued observation that the $L^{\infty}$ bound on the first order term of $\Delta_g$ can be relaxed to an $L^p$ bound so long as $p > n$. To explain the argument, suppose smooth functions $\phi_m$ are chosen so that $\phi_m \to \phi$ uniformly on $\partial \Omega$. Using Lemma \ref{thm:strong_solutions}, construct solutions $u_m \in W^{2,p}(\Omega)$ to $\Delta_g u_m = 0$, $u_m|_{\partial \Omega} = \phi_m$. Since $\phi_m \to \phi$ uniformly, the Alexandrov maximum principle \cite[Theorem 9.1]{MR1814364} implies $u_m$ is a uniformly Cauchy sequence in $\Omega$, so that there is a limit $u_m \to u \in C^0(\Omega)$. Interior regularity \cite[Theorem 9.11]{MR1814364} gives uniform boundedness in $W^{2,p}$ on compact subsets of $\Omega$, in turn giving $C^{1,\alpha'}$ convergence, $\alpha' < 1 - n/p$, from which we conclude that $\Delta_g u = 0$ on compact subsets of $\Omega$.
\end{proof}

We next aim to establish a harmonic atlas on $(M, g)$ and therefore optimal regularity of these coordinates, just as in \cite{MR644518}.
The existence of interior harmonic coordinates for $W^{1,p}$ metric tensors are well-known, see \cite[Chapter 3, Section 9]{MR1765330} for an exposition. We demonstrate the corresponding boundary coordinates here.

\begin{lemma}[Boundary harmonic coordinates for $W^{1,p}$ Riemannian manifolds]\label{thm:boundary_harmonic_existence}
	Suppose $(M, g)$ is a $W^{1,p}$ Riemannian manifold with boundary metric $h \in W^{1,p}(\partial M)$. To every point $z \in \partial M$ there exists a boundary harmonic coordinate chart centered at $z$.
\end{lemma}
\begin{proof}
	Suppose $U \ni z$ is a domain with a boundary chart $x$ centered at $z$, i.e. $x:U \to B_R^+$, $x = (x^1, \dots, x^n) = (x', x^n)$. Here $B_R^+$ is the intersection of the open ball of radius $R$ with $\{x^n > 0\}$. We apply the interior harmonic coordinate construction from \cite{MR1765330}. Thus, for small enough $r < R$ we construct $h$-harmonic coordinates $v$ on $\Sigma_r$, $v = (v^1, \dots, v^{n-1})$, such that $\Delta_h v^i = 0$ for $1 \le i \le n-1$ and $v \in W^{2,p}(\Sigma_r)$.

	Smooth the corners of $\overline{B^+_r}$ to obtain an approximating smooth domain $\Omega_r \subset \R^n_+$, with the additional requirement that $\partial \Omega_r \cap \{x^n = 0\} = \Sigma_r$. It is easy to do this, for instance, by using the graph of $x^n = \exp\left(-\frac{1}{(|x'|-r)^2}\right)$ on the region $\varepsilon + r > |x'| > r$, for some $\varepsilon > 0$, and then smoothly transitioning back to the standard hemisphere. This can be done in such a way that $\overline{\Omega}_r$ is uniformly close to $\overline{B^+_r}$, and such that $B^+_r \subset \Omega_r$ (think of a mushroom cap).

	We create an extension of each $v^i$ as follows. First, extend $v^i$ to a $W^{2,p}$ function on $\R^{n-1}$. Then, further extend this to a function $V^i \in W^{2,p}(\R^{n}_+)$. By the trace theorem, the restriction $V^i|_{\partial \Omega} \in W^{2-1/p, p}(\partial \Omega_r)$. See \cite[Theorems 2.7.2, 3.3.3, and 3.3.4]{MR3024598} for details on these extensions and traces, particularly Theorem 3.3.4 regarding an extension to $\R^{n-1}$, Theorem 2.7.2 regarding the extension from $\R^{n-1}$ to $\R^n$ (this is called a coretraction in their notation), and Theorem 3.3.3 for establishing the boundary trace to $\partial \Omega$.

	We now construct the full coordinates $u = (u^1, \dots, u^n)$ as solutions to appropriate Dirichlet problems on $\Omega_r$. To construct the tangential coordinates, we solve
	\begin{equation}\label{bndry_tangential}
		\begin{cases}
			\Delta_g u^i = 0 & \text{ in } \Omega_r,          \\
			u^i = V^i        & \text{ on } \partial \Omega_r,
		\end{cases}
	\end{equation}
	$1 \le i \le n-1$. To construct $u^n$, we solve
	\begin{equation}\label{bndry_normal}
		\begin{cases}
			\Delta_g u^n = 0 & \text{ in } \Omega_r,          \\
			u^n = x^n        & \text{ on } \partial \Omega_r.
		\end{cases}
	\end{equation}

	By Lemma \ref{thm:strong_solutions}, these Dirichlet problems admit unique strong solutions $u^1, \dots, u^n \in W^{2,p}(\Omega_r)$. Let us show the map $u = (u^1, \dots, u^n)$ is a diffeomorphism onto its image when restricted to a small enough half-ball about $z$. First, note that $u^n = 0$ on $\Sigma_r$ and $u^n = x^n > 0$ on $\partial \Omega_r \cap \{x^n > 0\}$. Applying the strong minimum principle for divergence-form equations (see for instance \cite[Theorem 8.19]
	{MR1814364}) we deduce that $u^n > 0$ in the interior of $\Omega_r$. Next, we use a version of the Hopf boundary point lemma due to Safonov \cite[Theorem 4.3 and Remark 4.4]{MR2667641} to conclude that the inner normal derivative is strictly positive at $z$. Regarding the tangential components $u^i$ ($1 \le i \le n-1$), the boundary condition $u^i(x', 0) = v^i(x')$ ensures that the tangential Jacobian matrix $\left( \frac{\partial u^i}{\partial x^j} \right)$ is non-singular at the origin, as $v$ forms a valid coordinate system on $\partial M$. Combining these facts, the full Jacobian matrix $Du$ is nonsingular, so that $u$ is a harmonic coordinate chart on a sufficiently small half-ball containing $z$.
\end{proof}

\subsection{Boundary regularity and the Neumann problem}
Of specific interest in our work is the following system of equations, due to Anderson and further studied in \cite{MR2096795}, for the metric and its inverse. Valid in harmonic coordinates, it includes Neumann conditions for the inverse metric components $g^{in}$, $1 \le i \le n$, and Dirichlet conditions for the tangential components $g_{\alpha\beta}$, $1 \le \alpha, \beta \le n-1$. Thus, on a smooth Riemannian manifold $(M,g)$, in a boundary harmonic coordinate system and writing $\partial_{\nu}$ for the inward normal vector field, one has

\begin{align}
	\Delta g_{\alpha\beta}    & = B_{\alpha\beta}(g, \partial g) - 2 \Ric(g)_{\alpha\beta}                       & \text{in } M, \label{eq:tangential_interior}          \\
	g_{\alpha\beta}           & = h_{\alpha\beta}                                                                & \text{on } \partial M, \label{eq:tangential_boundary} \\
	\Delta g^{in}             & = B^{in}(g, \partial g) - 2 \Ric(g)^{in}                                         & \text{in } M, \label{eq:normal_interior}              \\
	\partial_\nu g^{\alpha n} & = -(n-1)H g^{\alpha n} + \frac{1}{2\sqrt{g^{nn}}} g^{\alpha k} \partial_k g^{nn} & \text{on } \partial M, \label{eq:normal_boundary_1}   \\
	\partial_\nu g^{nn}       & = -2(n-1)H g^{nn}                                                                & \text{on } \partial M, \label{eq:normal_boundary_2}
\end{align}
where $B_{\alpha\beta}(g, \partial g)$ and $B^{in}(g, \partial g)$ are lower-order terms quadratic in $g$, $g^{-1}$, and $\partial g$.
We wish to study this system on a $W^{1,p}$ Riemannian manifold. Equations \eqref{eq:tangential_interior} - \eqref{eq:tangential_boundary} are straightforward to study with standard $L^p$ theory, e.g. with Lemmas~\ref{thm:local_boundary_w2p_nonhom} and~\ref{thm:strong_solutions}. Equations \eqref{eq:normal_interior} - \eqref{eq:normal_boundary_2} are more nuanced, as $\partial_{\nu} g^{in}$ is not a priori defined on $\partial M$ when $g \in W^{1,p}(M)$ (but see \cite[Section 2.2 and 2.3]{MR2096795}, and Proposition~\ref{thm:boundary_regularity} below). To begin, we recall the weak formulation of the Neumann problem and establish some a priori estimates.

\begin{definition}\label{def:weak_neumann}
	A function $u \in H^1(M)$ is a weak solution to the Neumann problem
	\begin{equation}
		\begin{cases}
			\Delta_g u = F      & \text{ in } M,          \\
			\partial_{\nu}u = G & \text{ on } \partial M,
		\end{cases}
	\end{equation}
	if, for every test function $\psi \in H^1(M)$ there holds
	\begin{equation} \label{eq:weak_neumann}
		\int_M \langle d u, d\psi \rangle_g \, dV_g = -(F, \psi)_{M} - (G, \psi)_{\partial M},
	\end{equation}
	where $(\cdot, \cdot)_M$ denotes the duality pairing between $H^{-1}(M)$ and $H^1(M)$, and $(\cdot, \cdot)_{\partial M}$ denotes the duality pairing between $H^{-1/2}(\partial M)$ and the trace space $H^{1/2}(\partial M)$.

	When the data are sufficiently regular, these dual pairings reduce to standard integrals over the manifold and its boundary,
	\begin{equation*}
		\int_M \langle d u, d\psi \rangle_g \, dV_g = -\int_M F \psi \, dV_g - \int_{\partial M} G \psi \, dS_g.
	\end{equation*}
\end{definition}

From this definition, we can derive $H^1$ a priori estimates.
\begin{lemma}[$H^1$ Neumann Estimates]\label{prop:neumann_w12_estimate}
	Let $(M,g)$ be a compact $W^{1,p}$ Riemannian manifold with boundary and $u \in H^1(M)$ be a weak solution to the Neumann problem
	\begin{equation}
		\begin{cases}
			\Delta_g u = F     & \text{ in } M,          \\
			\partial_\nu u = G & \text{ on } \partial M,
		\end{cases}
	\end{equation}
	with $F \in H^{-1}(M)$ and  $G \in H^{-1/2}(\partial M)$. Then there exists a constant $C > 0$, depending only on $n, M$, the ellipticity constants of $g$, and $\|g\|_{W^{1,p}(M)}$, such that
	\begin{equation} \label{eq:neumann_l2_estimate}
		\|u\|_{H^1(M)} \le C \left( \|F\|_{H^{-1}(M)} + \|G\|_{H^{-1/2}(\partial M)} + \|u\|_{L^2(M)} \right).
	\end{equation}
\end{lemma}

\begin{proof}
	Choosing the solution $u$ as the test function in the weak formulation of the Neumann problem, we obtain the energy identity
	\begin{equation} \label{eq:w12_energy}
		\int_M |d u|_g^2 \, dV_g = - (F, u)_M - (G, u)_{\partial M},
	\end{equation}
	whence
	\begin{equation} \label{eq:w12_lhs_bound}
		\|u\|_{H^1(M)}^2 = - (F, u)_M - (G, u)_{\partial M} + \|u\|_{L^2(M)}^2.
	\end{equation}
	On the right-hand side, we apply the Cauchy-Schwarz inequality and the trace theorem. The standard trace operator $T: H^1(M) \to H^{1/2}(\partial M)$ is bounded, so $\|u\|_{H^{1/2}(\partial M)} \le C_{tr} \|u\|_{H^1(M)}$. Thus:
	\begin{align}
		|(F, u)_M|            & \le \|F\|_{H^{-1}(M)} \|u\|_{H^1(M)},                                                                                \\
		|(G, u)_{\partial M}| & \le \|G\|_{H^{-1/2}(\partial M)} \|u\|_{H^{1/2}(\partial M)} \le C_{tr} \|G\|_{H^{-1/2}(\partial M)} \|u\|_{H^1(M)}.
	\end{align}
	Substituting these bounds into \eqref{eq:w12_lhs_bound} gives
	\begin{equation}
		\|u\|_{H^1(M)}^2 \le \left( \|F\|_{H^{-1}(M)} + C_{tr} \|G\|_{H^{-1/2}(\partial M)} \right) \|u\|_{H^1(M)} + \|u\|_{L^2(M)}^2.
	\end{equation}
	By applying Young's inequality ($ab \le \epsilon a^2 + \frac{1}{4\epsilon} b^2$) to the linear terms in $\|u\|_{H^1(M)}$ with $\epsilon = 1/2$, we absorb the $H^1$ norm into the left-hand side
	\begin{equation}
		\frac{1}{2} \|u\|_{H^1(M)}^2 \le C \left( \|F\|_{H^{-1}(M)}^2 + \|G\|_{H^{-1/2}(\partial M)}^2 \right) + \|u\|_{L^2(M)}^2,
	\end{equation}
	which implies the estimate.
\end{proof}

The weak Neumann problem can be localized as follows.
\begin{definition}[Local weak Neumann problem.]\label{def:local_weak_neumann}
	If $\mathcal O \subset M$ is an open set intersecting $\partial M$, write $U = \mathcal O \backslash \partial M$ and $\Sigma = \mathcal O \cap \partial M$. A function $u \in H^1(U)$ is a weak solution to the local weak Neumann problem
	\begin{equation}
		\begin{cases}
			\Delta_g u = F      & \text{ in } U,      \\
			\partial_{\nu}u = G & \text{ on } \Sigma,
		\end{cases}
	\end{equation}
	if
	\begin{equation} \label{eq:local_weak_neumann}
		\int_U \langle d u, d\psi \rangle_g \, dV_g = -(F, \psi)_{U} - (G, \psi)_{\Sigma}
	\end{equation}
	holds for each test function $\psi \in H^1_0(\mathcal O)$.
\end{definition}
We obtain a corresponding local $H^1$ estimate.
\begin{corollary}[Local $H^1$ Neumann Estimate]\label{cor:local_neumann_w12}
	Let $(M,g)$ be a $W^{1,p}$ Riemannian manifold with boundary. Suppose $u$ is a weak solution to the local Neumann problem as in Definition \ref{def:local_weak_neumann},
	with $F \in H^{-1}(U)$ and $G \in H^{-1/2}(\Sigma)$. Then for any compact subset $V \subset U\cup \Sigma$, there exists a constant $C > 0$, depending on $V, U$, the ellipticity constants of $g$, and $\|g\|_{W^{1,p}(U)}$, such that
	\begin{equation} \label{eq:local_neumann_w12_estimate}
		\|u\|_{H^1(V)} \le C \left( \|F\|_{H^{-1}(U)} + \|G\|_{H^{-1/2}(\Sigma)} + \|u\|_{L^2(U)} \right).
	\end{equation}
\end{corollary}

\begin{proof}
	Choose a smooth cutoff function $\eta \in C^\infty_c(U)$ such that $0 \le \eta \le 1$ and $\eta \equiv 1$ on $V$. We use the test function $\psi = \eta^2 u$. Since $\eta$ is compactly supported in $U \cup \Sigma$ and $u \in H^1(U)$, we have $\psi \in H^1_c(U)$, making it a valid test function for the local weak formulation.

	By the Sobolev product rule, $d\psi = \eta^2 du + 2\eta u \, d\eta$. Thus
	\begin{equation} \label{eq:caccioppoli_identity}
		\int_U \langle du, \eta^2 du + 2\eta u \, d\eta \rangle_g \, dV_g = - (F, \eta^2 u)_U - (G, \eta^2 u)_\Sigma.
	\end{equation}
	Expanding the left-hand side, we have
	\begin{equation*}
		\int_U \eta^2 |du|^2_g \, dV_g + 2 \int_U \eta u \langle du, d\eta \rangle_g \, dV_g.
	\end{equation*}
	Applying the Cauchy-Schwarz and Young's inequalities to the second term gives
	\begin{equation*}
		\left| 2 \int_U \eta u \langle du, d\eta \rangle_g \, dV_g \right| \le \frac{1}{2} \int_U \eta^2 |du|^2_g \, dV_g + 2 \int_U u^2 |d\eta|^2_g \, dV_g.
	\end{equation*}
	Together with \eqref{eq:caccioppoli_identity} we obtain
	\begin{equation*}
		\frac{1}{2} \int_U \eta^2 |du|^2_g \, dV_g - 2 \int_U u^2 |d\eta|^2_g \, dV_g \le \left| (F, \eta^2 u)_U \right| + \left| (G, \eta^2 u)_\Sigma \right|.
	\end{equation*}
	Next, we bound the data terms using the trace theorem, thus
	\begin{align*}
		\left| (F, \eta^2 u)_U \right|      & \le \|F\|_{H^{-1}(U)} \|\eta^2 u\|_{H^1(U)},                                                                                  \\
		\left| (G, \eta^2 u)_\Sigma \right| & \le \|G\|_{H^{-1/2}(\Sigma)} \|\eta^2 u\|_{H^{1/2}(\Sigma)} \le C_{\text{tr}} \|G\|_{H^{-1/2}(\Sigma)} \|\eta^2 u\|_{H^1(U)}.
	\end{align*}
	By the product rule and the triangle inequality, we have
	\begin{equation*}
		\|\eta^2 u\|_{H^1(U)} \le \|\eta du\|_{L^2(U)} + C_\eta \|u\|_{L^2(U)},
	\end{equation*}
	where $C_\eta$ depends on $\|\eta\|_{C^1}$. Applying Young's inequality with a sufficiently small $\epsilon > 0$ to the products on the right-hand side allows us to absorb the gradient term $\epsilon \|\eta du\|_{L^2(U)}^2$ into the left-hand side of the energy inequality. The remaining terms are bounded by $C(\|F\|_{H^{-1}(U)}^2 + \|G\|_{H^{-1/2}(\Sigma)}^2 + \|u\|_{L^2(U)}^2)$.
	Thus,
	\begin{equation*}
		\int_U \eta^2 |du|^2_g \, dV_g \le C \left( \|F\|_{H^{-1}(U)}^2 + \|G\|_{H^{-1/2}(\Sigma)}^2 + \|u\|_{L^2(U)}^2 \right).
	\end{equation*}
	Since $\eta \equiv 1$ on $V$, we have $\int_V |du|^2_g \, dV_g \le \int_U \eta^2 |du|^2_g \, dV_g$. Adding $\|u\|_{L^2(V)}^2 \le \|u\|_{L^2(U)}^2$ to both sides yields the desired local estimate \eqref{eq:local_neumann_w12_estimate}.
\end{proof}

We obtain further estimates to solutions of the weak Neumann problem by taking the viewpoint in \cite{MR2096795}, relying heavily on the layer potential theory developed by Mitrea and Taylor \cite{MR1765156, MR1846208, MR1953534}, which applies to subdomains of Riemannian manifolds whose boundaries are at least Lipschitz and with metric tensors that are at least H\"older continuous (other regularity scales are also studied as part of their program). Thus, we require a background Riemannian manifold $N$, smooth, compact, connected, and with a H\"older continuous metric tensor. On this background manifold, one considers the operator
\begin{equation}\label{eq:operator_L}
	L = \Delta - V,
\end{equation}
with $V \in L^{\infty}(N)$ and $V > 0$ on a set of positive measure somewhere on $N$. In our case, we take $N$ to be the smooth double of $M$, carrying with it a $W^{1,p}$ (thus $C^{0,\alpha}$) extension of the metric $g$, and we identify $M$ as a subset of $N$, taking $V = 0$ on $M$. The following result demonstrates the quantitative estimate associated to \cite[Proposition 5.5.1]{MR2096795}.

\begin{proposition}[$W^{1,q}$ Neumann Estimate]\label{prop:neumann_l_s_estimate}
	Suppose $(M,g)$ is a $W^{1,p}$ Riemannian manifold with boundary, $p > n$. Let $u \in H^1(M)$ be a weak solution to the Neumann problem
	\begin{equation}
		\begin{cases}
			\Delta_g u = F       & \text{ in } M,          \\
			\partial_{\nu} u = G & \text{ on } \partial M.
		\end{cases}
	\end{equation}
	Assume $F \in L^p(M)$ for $p > n$, and $G \in L^s(\partial M)$ for $1 < s < \infty$. Define $q = \frac{ns}{n-1}$. Then $u \in W^{1,q}(M)$, and the nontangential maximal function of the gradient, denoted $(\nabla u)^*$, belongs to $L^s(\partial M)$. Furthermore, there exists a constant $C > 0$, depending only on $n, p, s, M$, the ellipticity constants of $g$, and $\|g\|_{W^{1,p}(M)}$, such that
	\begin{equation} \label{eq:neumann_explicit_estimate}
		\|u\|_{W^{1,q}(M)} \le C \left( \|F\|_{L^p(M)} + \|G\|_{L^s(\partial M)} + \|u\|_{L^q(M)} \right),
	\end{equation}
	and the nontangential maximal function satisfies the quantitative bound
	\begin{equation} \label{eq:neumann_ntmf_estimate}
		\|(\nabla u)^*\|_{L^s(\partial M)} \le C \left( \|F\|_{L^p(M)} + \|G\|_{L^s(\partial M)} \right).
	\end{equation}
\end{proposition}

\begin{proof}
	First, extend $F$ by zero to $\tilde{F} \in L^p(N)$. Just as in \cite[Theorem 5.3.1]{MitreaTaylor2005}, for $\sigma < 1 - n/p$ we obtain $v \in C^{1,\sigma}(N)$ satisfying $L v = \tilde{F}$ with the bound
	\begin{equation}\label{eq:3.10_interior}
		\|v\|_{C^{1,\sigma}(N)} \le C\|\tilde{F}\|_{L^p(N)} = C\|F\|_{L^p(M)}.
	\end{equation}
	Because $V \equiv 0$ on $M$, $v$ satisfies $\Delta_g v = F$ on $M$.	The normal derivative $\partial_\nu v$ and the nontangential maximal function $(\nabla v)^*$ are  trivially in $L^\infty(\partial M) \subset L^s(\partial M)$, satisfying
	\begin{equation}\label{eq:3.10_v_bound}
		\|(\nabla v)^*\|_{L^s(\partial M)} + \|\partial_\nu v\|_{L^s(\partial M)} \le C\|F\|_{L^p(M)}.
	\end{equation}

	Now consider $w = u - v|_M$. By linearity, $w \in H^1(M)$ is a weak solution to the homogeneous equation with modified Neumann data
	\begin{equation}\label{eq:3.12}
		\begin{cases}
			\Delta_g w = 0                       & \text{in } M,          \\
			\partial_{\nu}w = G - \partial_\nu v & \text{on } \partial M.
		\end{cases}
	\end{equation}
	The modified boundary data $G - \partial_\nu v$ lies in $L^s(\partial M)$. From the Mitrea-Taylor layer potential theory (\cite[Section 4]{MR2130803} and \cite[Proposition 5.5.1]{MR2096795}), there exists a bounded Neumann solution operator $NI : L^s(\partial M) \to W^{1,q}(M)$. Setting $w_{NI} := NI(G - \partial_\nu v)$, the nontangential maximal function of the gradient is controlled by the $L^s$ boundary data, thus
	\begin{equation}\label{eq:3.13_ntmf}
		\|(\nabla w_{NI})^*\|_{L^s(\partial M)} \le C\|G - \partial_\nu v\|_{L^s(\partial M)} \le C\left(\|G\|_{L^s(\partial M)} + \|F\|_{L^p(M)}\right).
	\end{equation}
	By the boundedness of $NI$ mapping into $W^{1,q}(M)$, we analogously have
	\begin{equation}\label{eq:3.13}
		\|w_{NI}\|_{W^{1,q}(M)} \le C\|G - \partial_\nu v\|_{L^s(\partial M)} \le C\left(\|G\|_{L^s(\partial M)} + \|F\|_{L^p(M)}\right).
	\end{equation}

	Since $\Delta_g w = 0$, $w$ differs from $w_{NI}$ by an element $\eta \in \ker(\Delta_g)$; thus $\eta$ is a constant. Since $u = v + w_{NI} + \eta$, we get
	\begin{align}
		\|(\nabla u)^*\|_{L^s(\partial M)} & \le \|(\nabla v)^*\|_{L^s(\partial M)} + \|(\nabla w_{NI})^*\|_{L^s(\partial M)} + 0 \nonumber \\
		                                   & \le C\|F\|_{L^p(M)} + C\left(\|G\|_{L^s(\partial M)} + \|F\|_{L^p(M)}\right) \nonumber         \\
		                                   & \le C\left(\|F\|_{L^p(M)} + \|G\|_{L^s(\partial M)}\right), \label{eq:3.15_ntmf}
	\end{align}
	which establishes \eqref{eq:neumann_ntmf_estimate}.

	To establish the $W^{1,q}$ bound, we isolate $\eta = w - w_{NI} = u - v - w_{NI}$. Since $\eta$ is constant, we see that $\|\eta\|_{W^{1,q}(M)} = \|\eta\|_{L^q(M)}$. Then we get
	\begin{equation}\label{eq:3.14}
		\|\eta\|_{W^{1,q}(M)} \le \|u\|_{L^q(M)} + \|v\|_{L^q(M)} + \|w_{NI}\|_{L^q(M)}.
	\end{equation}
	Finally, we arrive at
	\begin{align*}
		\|u\|_{W^{1,q}(M)} & \le \|v\|_{W^{1,q}(M)} + \|w_{NI}\|_{W^{1,q}(M)} + \|\eta\|_{W^{1,q}(M)}                                                                  \\
		                   & \le C\|v\|_{C^{1,\sigma}(M)} + \|w_{NI}\|_{W^{1,q}(M)} + C\left(\|u\|_{L^q(M)} + \|v\|_{C^{1,\sigma}(M)} + \|w_{NI}\|_{W^{1,q}(M)}\right) \\
		                   & \le C\left(\|F\|_{L^p(M)} + \|G\|_{L^s(\partial M)} + \|u\|_{L^q(M)}\right),
	\end{align*}
	where terms involving $v$ and $w_{NI}$ have been bounded by the problem data using \eqref{eq:3.10_interior} and \eqref{eq:3.13}. This proves \eqref{eq:neumann_explicit_estimate}.
\end{proof}

It follows from \cite{MR1953534} that the non-tangential maximal function bound in the previous result implies that the gradient of $u$ has an almost-everywhere trace that
is bounded in $L^s$. We record the estimate below.

\begin{corollary}[Boundary Trace of the Gradient]\label{cor:neumann_gradient_trace}
	Under the assumptions of Proposition \ref{prop:neumann_l_s_estimate}, the gradient of the weak solution, $\nabla u$, admits a nontangential limit almost everywhere on $\partial M$.
	Denoting this boundary trace by $\nabla u|_{\partial M}$, we have $\nabla u|_{\partial M} \in L^s(\partial M)$, satisfying
	\begin{equation} \label{eq:neumann_gradient_trace_estimate}
		\|\nabla u|_{\partial M}\|_{L^s(\partial M)} \le C \left( \|F\|_{L^p(M)} + \|G\|_{L^s(\partial M)} \right),
	\end{equation}
	where $C > 0$ depends only on $n, p, s, M$, the ellipticity constants of $g$, and $\|g\|_{W^{1,p}(M)}$.
\end{corollary}

\begin{proof}
	We utilize the decomposition $u = v + w_{NI} + \eta$ from the proof of Proposition \ref{prop:neumann_l_s_estimate}. Since $v \in C^{1,\sigma}(N)$ and $\eta$ is constant, we only need to study $w_{NI}$ carefully. We have that $w_{NI} = NI(G - \partial_\nu v)$ is the solution constructed via layer potential operators for the homogeneous
	equation $\Delta_g w = 0$ with modified Neumann data $G - \partial_\nu v \in L^s(\partial M)$. From Proposition \ref{prop:neumann_l_s_estimate}, $(\nabla w_{NI})^*$ is bounded in $L^s$, which implies it is finite almost everywhere. Thus, by \cite[Theorem 4.2]{MR1953534}, $\nabla w_{NI}|_{\partial M}$ exists pointwise almost everywhere on $\partial M$ (see also \cite[Proposition 3.1]{MR1846208} for a more detailed argument assuming a Lipschitz bound on the metric. The Lipschitz bound is removed in \cite{MR1953534} while the Fatou-type argument in \cite{MR1846208} stays essentially the same.)

	By linearity, the nontangential limit of $\nabla u$ exists almost everywhere on $\partial M$:
	\begin{equation}
		\nabla u|_{\partial M}(x) = \lim_{\substack{y \to x \\ y \in \Gamma_x}} \nabla u(y) \quad \text{for a.e. } x \in \partial M,
	\end{equation}
	where $\Gamma_x = \{y \in M \setminus \partial M : d(y,x) \le 2\, d(y, \partial M)\}$ is a nontangential approach cone at $x$. To establish the quantitative bound, we observe that the nontangential limit is bounded by the maximal function, so that
	\begin{equation}
		|\nabla u|_{\partial M}(x)| \le (\nabla u)^*(x) \quad \text{for a.e. } x \in \partial M.
	\end{equation}
	Taking the $L^s(\partial M)$ norm of both sides and applying the quantitative bound \eqref{eq:neumann_ntmf_estimate} from Proposition \ref{prop:neumann_l_s_estimate} gives
	\begin{equation}
		\|\nabla u|_{\partial M}\|_{L^s(\partial M)} \le \|(\nabla u)^*\|_{L^s(\partial M)} \le C \left( \|F\|_{L^p(M)} + \|G\|_{L^s(\partial M)} \right).
	\end{equation}
\end{proof}

\begin{corollary}[Local $W^{1,q}$ Neumann Estimate]\label{prop:local_neumann_l_s}
	Suppose $g$ is a $W^{1,p}$ Riemannian metric on $B_R^+$. Let $u \in W^{1,p}(B_R^+)$ be a weak solution to the local Neumann problem
	\begin{equation}
		\begin{cases}
			\Delta_g u = F       & \text{ in } B_R^+,    \\
			\partial_{\nu} u = G & \text{ on } \Sigma_R.
		\end{cases}
	\end{equation}
	Assume $F \in L^p(B_R^+)$ and $G \in L^s(\Sigma_R)$ for $1 < s < \infty$. Define $q = \frac{ns}{n-1}$. Then for any $r < R$, $u \in W^{1,q}(B_r^+)$, and the nontangential maximal function of its gradient, $(\nabla u)^*$, belongs to $L^s(\Sigma_r)$. Furthermore, there exists a constant $C > 0$, depending strictly on $n, p, s, r, R$, the ellipticity constants of $g$, and $\|g\|_{W^{1,p}(B_R^+)}$, such that:
	\begin{equation} \label{eq:local_neumann_explicit_estimate}
		\|u\|_{W^{1,q}(B_r^+)} + \|(\nabla u)^*\|_{L^s(\Sigma_r)} \le C \left( \|F\|_{L^p(B_R^+)} + \|G\|_{L^s(\Sigma_R)} + \|u\|_{W^{1,p}(B_R^+)} \right).
	\end{equation}
\end{corollary}

\begin{proof}
	If $r < R$, construct a smooth domain $\Omega_{r}$ (parallel to Lemma \ref{thm:boundary_harmonic_existence}) so that $B_r^+ \subset \Omega_r \subset  B_R^+$ and with the property that $\partial \Omega_r \cap \{x^n = 0\} = \Sigma_r$. We write $N$ for the smooth double of $\Omega_r$ and extend $g$ to $W^{1,p}(N)$, viewing $\Omega_r \subset N$ as a smooth subdomain. Next, choose a smooth cutoff function $\eta \in C^\infty(N)$ such that $0 \le \eta \le 1$, $\eta \equiv 1$ on $B_{\rho}(0) \subset N$, $\rho < r$, and $\operatorname{supp}(\eta) \subset B_{\rho'}(0)$ for $\rho < \rho' < r$. Define $v = \eta u$, extended by zero to all of $N$. For any test function $\psi \in H^1(\Omega_r)$, we note that $v$ satisfies the global weak formulation
	\begin{equation*}
		\int_{\Omega_r} \langle dv, d\psi \rangle_g \,dV_g = -\int_{\Omega_r} \tilde{F} \psi \,dV_g - \int_{\partial \Omega_r} \tilde{G} \psi \,dS_g,
	\end{equation*}
	where the data are given by
	\begin{align*}
		\tilde{F} & = \eta F + 2\langle du, d\eta \rangle_g + u\Delta_g \eta, \\
		\tilde{G} & =
		\begin{cases}
			\eta G + u\partial_\nu \eta & \text{on } \Sigma_{r},                               \\
			0                           & \text{on } \partial \Omega_{r} \setminus \Sigma_{r}.
		\end{cases}
	\end{align*}

	The hypotheses on $g$, $u$, and $\eta$ imply that $\tilde F \in L^p(\Omega_r)$ and $\tilde G \in L^s(\partial \Omega_r)$. Specifically, the coefficients of $\Delta_g$ are at least $L^p(\Omega_{r})$ and both $u$ and $g$ are $C^{0,\alpha}(\Omega_r)$ for $\alpha = 1 - n/p$, as well as $W^{1,p}(\Omega_r)$, giving $\tilde F \in L^p(\Omega_r)$. Similarly, the condition on $g$ implies $\partial_{\nu} \in C^{0,\alpha}(\partial \Omega_r)$ allowing us to conclude $\tilde G \in L^s(\partial \Omega_r)$. Applying Proposition \ref{prop:neumann_l_s_estimate} to $v$ on $\Omega_r$ implies $v \in W^{1,q}(\Omega_r)$, $(\nabla v)^* \in L^s(\partial \Omega_r)$, and:
	\begin{equation*}
		\|v\|_{W^{1,q}(\Omega_r)} + \|(\nabla v)^*\|_{L^s(\partial \Omega_r)} \le C \left( \|\tilde{F}\|_{L^p(\Omega_r)} + \|\tilde{G}\|_{L^s(\partial \Omega_r)} + \|v\|_{L^q(\Omega_r)} \right).
	\end{equation*}
	Using the bounds
	\begin{align*}
		\|\tilde{G}\|_{L^s(\partial \Omega_r)} & \le C\left(\|G\|_{L^s(\Sigma_R)} + \|u\|_{W^{1,p}(B_R^+)}\right) \quad \text{and} \\
		\|\tilde{F}\|_{L^p(\Omega_r)}          & \le C \left(\|F\|_{L^p(B_R^+)} + \|u\|_{W^{1,p}(B_R^+)}\right),
	\end{align*}
	and bounding $\|v\|_{L^q(\Omega_r)}$ via the continuous embeddings $W^{1,p} \hookrightarrow L^\infty \hookrightarrow L^q$, we arrive at
	\begin{equation*}
		\|v\|_{W^{1,q}(\Omega_r)} + \|(\nabla v)^*\|_{L^s(\partial \Omega_r)} \le C \left( \|F\|_{L^p(B_R^+)} + \|G\|_{L^s(\Sigma_R)} + \|u\|_{W^{1,p}(B_R^+)} \right).
	\end{equation*}
	Since $v \equiv u$ on $B_{\rho}^+$, and $\rho$ and $r$ were arbitrary, the stated local estimate follows.
\end{proof}

We conclude this section by establishing boundary regularity for $W^{1,p}$ Riemannian manifolds satisfying $L^{\infty}$ curvature bounds. This is analogous to \cite[Theorem 2.1]{MR2096795}, which shows $C^2_*$ regularity under the stronger assumption that $H$ is a Lipschitz function on $\partial M$.

\begin{proposition}[Regularity for $W^{1,p}$ Riemannian Manifolds]\label{thm:boundary_regularity}
	Let $(M, g)$ be a $W^{1,p}$ Riemannian manifold with boundary, $p > n$, satisfying
	\[ \|\Ric(g)\| \le K, \quad \|\Ric(h)\| \leq K, \quad \text{and} \quad \|H\| \leq H_0. \]
	Then $(M,g)$ is a $W^{1,q}$ Riemannian manifold for any $q > n$. Moreover, the gradient of $g$ admits a pointwise a.e.\ trace on $\partial M$ belonging to $L^q(\partial M)$ for any $q < \infty$. Specifically, in any boundary harmonic coordinate atlas, the components $g^{in}$ satisfy the Neumann problem \eqref{eq:normal_interior}--\eqref{eq:normal_boundary_2} in the weak sense given by \eqref{eq:weak_neumann} with
	\begin{equation} \label{eq:source_F}
		F := F^{in} = B^{in}(g, \partial g) - 2 \Ric(g)^{in},
	\end{equation}
	and boundary data
	\begin{equation} \label{eq:source_G}
		G := G^{in} =
		\begin{cases}
			-(n-1)H g^{\alpha n} + \frac{1}{2\sqrt{g^{nn}}} g^{\alpha k} \partial_k g^{nn} & \text{if } i = \alpha < n, \\
			-2(n-1)H g^{nn}                                                                & \text{if } i=n.
		\end{cases}
	\end{equation}
\end{proposition}

\begin{proof}
	The proof adapts the arguments of \cite[Theorem 2.1 and Lemmas 2.2.3--2.3.2]{MR2096795} to our relaxed $L^\infty$ mean curvature bound, utilizing the local $W^{1,q}$ Neumann estimates established above. First, the assumption $\|\Ric(g)\| \le K$ combined with standard interior elliptic regularity implies $g \in W^{2,q}_{\mathrm{loc}}(M \setminus \partial M)$ for any $q < \infty$. It remains to establish regularity at the boundary.

	Let $z \in \partial M$. By Lemma \ref{thm:boundary_harmonic_existence}, we construct a boundary harmonic coordinate chart $x: U \to B_R^+ \subset \mathbb{R}^n_+$ centered at $z$.
	If $p < 2n$, we bootstrap the regularity on a sequence of shrinking half-balls. The terms $F^{in}$ defined in \eqref{eq:source_F} satisfy $F^{in} \in L^{p/2}(B_R^+)$. We first analyze the tangential components $g_{\alpha\beta}$ ($1 \leq \alpha, \beta \leq  n-1$). Defining the corresponding tangential term $F_{\alpha\beta} := B_{\alpha\beta}(g, \partial g) - 2 \Ric(g)_{\alpha\beta} \in L^{p/2}(B_R^+)$, and noting that $\|\Ric(h)\| \leq K$ implies the boundary metric $h$ is in $W^{2,q}(\partial M)$ for any $q < \infty$, the components $g_{\alpha\beta}$ satisfy a local Dirichlet problem
	\begin{equation}
		\begin{cases}
			\Delta_g g_{\alpha\beta} = F_{\alpha\beta} & \text{ in } B_R^+,    \\
			g_{\alpha\beta} = h_{\alpha\beta}          & \text{ on } \Sigma_R.
		\end{cases}
	\end{equation}
	Local elliptic regularity gives $g_{\alpha\beta} \in W^{1, s_1}(B_{R_1}^+)$ for a slightly smaller radius $R_1 < R$, where the improved integrability exponent $s_1$ satisfies $\frac{1}{s_1} = \frac{1}{p/2} - \frac{1}{n}$. (If $p = 2n$, standard embedding allows us to bypass this equation and select $s_1$ arbitrarily large).

	Next, we verify the weak formulation \eqref{eq:local_weak_neumann} for the normal component $g^{nn}$ on $B_{R_1}^+$. Following the level-set argument in \cite[Lemma 2.2.3]{MR2096795}, the Neumann data on the interior level sets $\{ x \mid x^n = c\}$ are well-defined and pass to the limit as $c \to 0$ in $W^{-1/p,p}(\Sigma_{R_1})$. Because the mean curvature satisfies $H \in L^\infty(\Sigma_{R_1})$, we obtain a local weak Neumann problem on $B_{R_1}^+$ with $F^{nn} \in L^{p/2}(B_{R_1}^+)$ and boundary data $G^{nn} = -2(n-1)H g^{nn} \in L^{\infty}(\Sigma_{R_1})$.

	Applying local Neumann estimates (cf.\ \cite[Lemma 2.3.3]{MR2096795}) to this weak formulation, we conclude that $g^{nn} \in W^{1, s_1-\epsilon}(B_{R_2}^+)$ for any $\epsilon > 0$ and $R_2 < R_1$, and that its gradient admits a boundary trace $\nabla g^{nn}|_{\Sigma} \in L^{s_1-\epsilon}(\Sigma_{R_2})$.

	This gradient trace allows us to process the mixed components $g^{\alpha n}$. Because $\nabla g^{nn}|_{\Sigma} \in L^{s_1-\epsilon}(\Sigma_{R_2})$, the boundary data defined in \eqref{eq:source_G} satisfies $G^{\alpha n} \in L^{s_1-\epsilon}(\Sigma_{R_2})$. Satisfying the local weak Neumann problem with these bounds yields $g^{\alpha n} \in W^{1, s_1-\epsilon}(B_{R_3}^+)$.

	Together, these steps demonstrate that the full metric tensor satisfies $g_{ij} \in W^{1, s_1-\epsilon}(B_{R_3}^+)$. Since $s_1 - \epsilon > p$, the integrability strictly improves. Iterating this bootstrap argument finitely many times on successively shrinking half-balls, we obtain a radius $\rho > 0$ such that $g_{ij} \in W^{1,p'}(B_{\rho}^+)$ for an exponent $p' > 2n$.

	Once $p' > 2n$, the interior terms satisfy $F^{in} \in L^{p'/2}(B_{\rho}^+)$ with $p'/2 > n$. We may now apply the strong local estimates. For the normal component, the boundary data remains $G^{nn} \in L^\infty(\Sigma_{\rho})$. Applying the Local $W^{1,q}$ Neumann Estimate (Corollary \ref{prop:local_neumann_l_s}) to $g^{nn}$ guarantees that for any $\rho_1 < \rho$ and any $q < \infty$, $g^{nn} \in W^{1,q}(B_{\rho_1}^+)$, and the nontangential limit of its gradient belongs to $L^q(\Sigma_{\rho_1})$.

	With $\nabla g^{nn}|_\Sigma \in L^q(\Sigma_{\rho_1})$ for any $q < \infty$, the mixed boundary data now satisfies $G^{\alpha n} \in L^q(\Sigma_{\rho_1})$. A second application of the Local $W^{1,q}$ Neumann Estimate to $g^{\alpha n}$ immediately yields $g^{\alpha n} \in W^{1,q}(B_{\rho_2}^+)$ for any $\rho_2 < \rho_1$, accompanied by the corresponding gradient trace estimate in $L^q(\Sigma_{\rho_2})$. Since $z \in \partial M$ was arbitrary, these local bounds patch together to confirm $(M,g)$ is globally a $W^{1,q}$ Riemannian manifold for any $q > n$.
\end{proof}

\subsection{Convergence Theory}
We begin with a few definitions of convergence used in this paper.

\begin{definition}[$W^{k,p}$ Convergence]\label{def:wkp_convergence}
	A sequence $(M_i,g_i)$ of compact Riemannian manifolds with boundary converges to a limit $W^{k,p}$ Riemannian manifold $(M,g)$ in the $W^{k,p}$ topology, $k > n/p$, if for large $i$ there exists a sequence of $W^{k+1, p}$ diffeomorphisms $F_i: M \to M_i$ such that the pull-back metrics $F_i^*g_i$ converge to $g$ in $W^{k,p}(M)$.
\end{definition}

\begin{definition}[Local Convergence]\label{def:local_convergence}
	A sequence of Riemannian manifolds with boundary $(M_i, g_i)$ converges to a limit manifold $(M, g)$ in the $W^{k,p}_{\mathrm{loc}}(M \setminus \partial M)$ topology, for $k > n/p$, if for every compact set $V \subset M \setminus \partial M$, there exists a relatively compact open set $U$ such that $V \subset U \subset \overline{U} \subset M \setminus \partial M$, and for sufficiently large $i$, there exist $W^{k+1,p}$ diffeomorphisms onto their image
	\[
		F_i: U \to F_i(U) \subset M_i \setminus \partial M_i
	\]
	such that the pull-back metrics $F_i^*g_i$ converge to $g$ in the $W^{k,p}(V)$ topology.
\end{definition}

The blowup argument used in Theorem \ref{precompact} involves a sequence of compact manifolds converging to a noncompact limit. In that situation, as in \cite{MR2096795}, we use a notion of pointed convergence analogous to the local convergence defined above.

\begin{definition}[Pointed Convergence]\label{def:pointed_convergence}
	A sequence of pointed Riemannian manifolds with boundary $(M_i, g_i, p_i)$ converges to a pointed limit manifold $(M, g, p)$ in the pointed $W^{k,p}$ topology, for $k > n/p$, if there exist sequences of real numbers $r_i \to \infty$ and $s_i > r_i$, along with open sets $U_i \subset M_i$ and $V_i \subset M$, such that the following conditions hold:
	\begin{enumerate}
		\item The open sets capture arbitrarily large metric balls around the base points:
		      \begin{align*}
			      B_{r_i}(p_i) & \subset U_i \subset B_{s_i}(p_i), & B_{r_i}(p) & \subset V_i \subset B_{s_i}(p).
		      \end{align*}
		\item There exist $W^{k+1,p}$ diffeomorphisms $F_i: V_i \to U_i$. (Note that $F_i(V_i \cap \partial M) = U_i \cap \partial M_i$ is automatic).
		\item For every precompact open set $\Omega \subset M$, for all sufficiently large $i$, we have $\Omega \subset V_i$, and the pull-back metrics $F_i^*g_i$ converge to $g$ in the $W^{k,p}(\Omega)$ topology.
		\item The base points converge: $F_i^{-1}(p_i) \to p$ in $M$.
	\end{enumerate}
\end{definition}

Similar definitions hold for convergence in other function spaces. For an introduction to the convergence theory of Riemannian manifolds, see for instance \cite{gromov2006metric, MR2243772, MR2684779, MR1620864, Sormani2012}. We recall the definition of the harmonic radius, introduced by Anderson in \cite{MR1074481}.

\begin{definition}[$W^{k,p}$ Harmonic Radius]\label{def:harmonic_radius}
	Suppose $M$ is a $W^{k,p}$ Riemannian manifold. Fix a point $x_0 \in M \setminus \partial M$ and a constant $Q > 1$. Define the $W^{k,p}$-harmonic radius at $x_0$, denoted $r^Q_h(x_0) := r_h(x_0)$, as the supremum over all radii $r < \dist(x_0, \partial M)$ such that there exists a harmonic coordinate system $x = (x^1, \dots, x^n)$ centered at $x_0$. We require this chart to contain the geodesic ball $B_r(x_0)$, and the metric components $g_{ij}$ in these coordinates to satisfy
	\begin{align}
		Q^{-1}\delta_{ij} \leq g_{ij}                                & \leq Q\delta_{ij}, \label{harmonic} \\
		r^{|\sigma|-n/p} \|\partial^{\sigma}g_{ij}\|_{L^p(B_r(x_0))} & \leq Q, \label{coord}
	\end{align}
	for each multi-index $\sigma$ with $1 \leq |\sigma| \leq k$. Here, $\delta_{ij}$ denotes the standard Euclidean metric.
\end{definition}

The conditions \eqref{harmonic} and \eqref{coord} are invariant under simultaneous rescaling of the metric and the coordinates, meaning that $r_h$ scales as a distance function: $r_h(\lambda^2 g) = \lambda r_h(g)$. We'll introduce here a corresponding definition of a $W^{k,p}$ boundary harmonic radius.

\begin{definition}[$W^{k,p}$ Boundary Harmonic Radius]\label{def:boundary_harmonic_radius}
	Let $(M,g)$ be a $W^{k,p}$ Riemannian $n$-manifold with boundary. Suppose the boundary metric $h \in W^{k,p}(\partial M)$ and fix $Q > 1$. For a boundary point $z \in \partial M$, we define the $W^{k,p}$-boundary harmonic radius at $z$, denoted $r^Q_{bh}(z) := r_{bh}(z)$, as the supremum over all radii $r > 0$ such that there exists a boundary harmonic coordinate chart $x = (x^1, \dots, x^n)$ defined on the ball $B_r(z) \subset M$ satisfying the following properties:

	The components of the metric $g_{ij}$ satisfy
	\begin{align}
		Q^{-1}\delta_{ij} \leq g_{ij}                             & \leq Q\delta_{ij}, \quad \label{eq:interior_C0_bound} \\
		r^{|\sigma|-n/p} \|\partial^\sigma g_{ij}\|_{L^p(B_r(z))} & \leq Q, \label{eq:interior_Wkp_bound}
	\end{align}
	on $B_r(z)$ for each multi-index $\sigma$ with $1 \leq |\sigma| \leq k$.

	The components of the induced boundary metric $h_{\alpha\beta}$ satisfy
	\begin{align}
		Q^{-1}\delta_{\alpha\beta} \leq h_{\alpha\beta}                                        & \leq Q\delta_{\alpha\beta}, \label{eq:boundary_C0_bound} \\
		r^{|\sigma|-(n-1)/p} \|\partial^\sigma h_{\alpha\beta}\|_{L^p(\partial M \cap B_r(z))} & \leq Q, \label{eq:boundary_Wkp_bound}
	\end{align}
	on $\partial M \cap B_r(z)$ for each multi-index $\sigma$ with $1 \leq |\sigma| \leq k$. As in the interior case, the global boundary harmonic radius is defined as $r_{bh} := \inf_{z \in \partial M} r_{bh}(z)$.
\end{definition}

\begin{remark}\label{rmk:trace_regularity}
	If a metric $g$ belongs to $W^{k,p}(M)$, the induced boundary metric $h$ belongs to the fractional Sobolev space $W^{k-1/p, p}(\partial M)$. Consequently, a sharper definition of the boundary harmonic radius for a general $W^{k,p}$ metric would require replacing the $W^{k,p}$ bounds in \eqref{eq:boundary_Wkp_bound} with the corresponding $W^{k-1/p, p}$ Gagliardo seminorm. It is possible to do this, e.g.\ to define $W^{s,p}$, $s \in \R$, $s > n/p$, harmonic radii and to discuss their continuity, which would also involve bringing in $W^{s,p}$ elliptic estimates. However, this is not necessary for the classes of manifolds considered in this paper. See also Remark \ref{rmk:stratified_convergence}.
\end{remark}

Lemma \eqref{thm:boundary_harmonic_existence} shows that the $W^{1,p}$ global boundary harmonic radius $r_{bh} > 0$ on a compact $W^{1,p}$ Riemannian manifold with boundary. It is known that the interior $W^{1,p}$ harmonic radius is continuous with respect to strong $W^{1,p}$ convergence, see \cite{MR1158336} for an exposition. We show here that the boundary harmonic radius also has this property.

\begin{proposition} \label{thm:harmonic_radius_continuity}
	Let $M$ be a compact, smooth $n$-manifold with  boundary $\partial M$, and fix $p > n$ and $Q > 1$. Suppose $\{g_m\}_{m=1}^\infty$ is a sequence of Riemannian metrics on $M$ converging strongly in $W^{1,p}(M)$ to a metric $g$. Assume further that the induced boundary metrics $h_m$ converge strongly in $W^{1,p}(\partial M)$ to the limit boundary metric $h$. Then, for any $z \in \partial M$, the $W^{1,p}$ boundary harmonic radius converges,
	\begin{equation}
		\lim_{m \to \infty} r_{bh}(z, g_m) = r_{bh}(z, g).
	\end{equation}
	Furthermore, the global boundary harmonic radius converges,
	\begin{equation}
		\lim_{m \to \infty} r_{bh}(M, g_m) = r_{bh}(M, g).
	\end{equation}
\end{proposition}

\begin{proof}
	Write $r^+ = \limsup_{m \to \infty} r_{bh}(z, g_m)$ and assume $r^+ > 0$. First, we show that $r^+ \leq r_{bh}(z, g)$. Choose a subsequence (still denoted $g_m$) such that $r_{bh}(z, g_m) \to r^+$, and consider boundary harmonic coordinate charts $u_m: B_{r_m}(z) \to \mathbb{R}^n_+$ that satisfy the bounds in Definition \ref{def:boundary_harmonic_radius}. Write $\Sigma_{r_m} := B_{r_m}(z) \cap \partial M$.

	For any $r < r^+$, we cover $B_{r}(z)$ by finitely many boundary or interior charts. Applying Lemma~\ref{thm:local_boundary_w2p_nonhom} or the corresponding interior estimate in these charts and summing over the cover, we get, for large enough $m$ and $r < r' < r^+$,
	\begin{equation}\label{eq:xm_elliptic_estimate}
		\|u^i_m\|_{W^{2,p}(B_{r})} \le C \left( \|u^i_m\|_{W^{2-1/p,p}(\Sigma_{r'})} + \|u^i_m\|_{L^p(B_{r'})}\right),
	\end{equation}
	where the norms are taken on $(M,g)$ and $C$ only depends on $n, p, r, r'$, the ellipticity constants of $g$, and $\|g\|_{W^{1,p}(M)}$. Definition \ref{def:boundary_harmonic_radius} implies that $\|u_m\|_{L^p(B_{r'})}$ is uniformly bounded. An estimate similar to \eqref{eq:xm_elliptic_estimate} on $\partial M$, noting that the tangential coordinates are harmonic with respect to $h_m$, shows that the $W^{2-1/p,p}(\Sigma_{r'})$ norms of the tangential coordinates are also uniformly bounded. (In fact, since the harmonic radius on $\partial M$ is continuous, we get that $v_m = u_m |_{\Sigma_{r'}}$ converge in $W^{2,p}(\Sigma_{r'})$ to coordinates $v$ that are harmonic with respect to $h$.)

	Since $r'$ is arbitrary, and noting that $u_m^n = 0$ on $\partial M$, we get that $u^i_m$ is uniformly bounded in $W^{2,p}(B_r)$ for any $r < r^+$. After passing to a subsequence we conclude there is a limit $u:B_r \to \mathbb{R}^n_+$, with $u_m \to u$ weakly in $W^{2,p}(B_r)$ and strongly in $C^{1,\alpha'}(B_r)$ for $\alpha' < 1 - n/p$.  Next, applying the same local elliptic estimates as before to the difference $u_m^i - u^i$, we obtain
	\begin{equation}\label{eq:xm_x_diff_elliptic_estimate}
		\begin{split}
			\|u^i_m - u^i\|_{W^{2,p}(B_r)} & \le C \Big( \| (\Delta_{g_m} - \Delta_g)u^i \|_{L^p(B_{r'})}                               \\
			                               & \qquad + \|u^i_m - u^i\|_{W^{2-1/p,p}(\Sigma_{r'})} + \|u^i_m - u^i\|_{L^p(B_{r'})} \Big).
		\end{split}
	\end{equation}
	The first term on the right-hand side goes to zero because $g_m \to g$ strongly in $W^{1,p}$, the second term goes to zero due to the convergence of the harmonic coordinates on the intrinsic boundary, and the third term goes to zero due to the compact embedding $W^{2,p} \hookrightarrow L^p$. Thus, we conclude $u_m \to u$ strongly in $W^{2,p}(B_r)$ for any $r < r^+$. Together with the strong convergence of $g_m \to g$ and $h_m \to h$ in $W^{1,p}(M)$, this implies each function $u^i$ is $g$-harmonic on the interior and the tangential components are $h$-harmonic on the boundary. The conditions of Definition \ref{def:boundary_harmonic_radius} are preserved under this limit, and since the convergence is stronger than $C^1$, Definition \ref{def:boundary_harmonic_radius} implies that $u$ is a valid coordinate chart on $M$. Thus $r_{bh}(z,g) \geq \limsup_{m \to \infty} r_{bh}(z, g_m)$.

	Next, we show that $r^- = \liminf_{m \to \infty} r_{bh}(z, g_m) \geq r_{bh}(z, g)$. Let $r < r_{bh}(z, g)$. Choose a radius $r'$ such that $r < r' < r_{bh}(z, g)$. Then there exists a boundary harmonic coordinate chart $u: U  \to B_R^+$, $B_{r'}(z) \Subset U$ and with $B^+_R$ the Euclidean upper half-ball of radius $R$, satisfying the interior $Q$-bounds \eqref{eq:interior_C0_bound}--\eqref{eq:interior_Wkp_bound} and the boundary $Q$-bounds \eqref{eq:boundary_C0_bound}--\eqref{eq:boundary_Wkp_bound} on $B_{r'}(z)$.

	Choose $R' < R$ so that $B_{r'}(z) \Subset u^{-1}(B^+_{R'})$. As in Lemma~\ref{thm:boundary_harmonic_existence}, we consider a smoothed domain $\Omega_{R'}$ satisfying $B_{R'}^+ \subset \Omega_{R'} \subset B_R^+$ and with $\Sigma_{R'} := \partial \Omega_{R'}\cap \overline B^+_{R} = \overline B_{R'}^+ \cap \{x^n = 0\}$. To construct boundary harmonic coordinates $u_m$ for $g_m$ on $B_r(z)$, we solve the Dirichlet problem on
	$\Omega_{R'}$. First, construct $h_m$-harmonic coordinates $(v^1_m, \ldots, v^{n-1}_m)$ in $W^{2,p}(\Sigma_{R'})$ with boundary values $v^i = u^i$ on $\Gamma_{R'} := \partial \Sigma_{R'} \subset \R^{n-1}$, with the property that $\|v_m - u\|_{W^{2,p}(\Sigma_{R'})} \to 0$. This is possible by the continuity of the harmonic radius on manifolds without boundary (although the argument is essentially the same as the one to follow). Define the piecewise boundary data $V_m^i$ on $\partial \Omega_{R'}$ by setting $V_m^i = v_m^i$ on $\Sigma_{R'}$ and $V_m^i = u^i$ on the curved portion $\partial \Omega_{R'} \setminus \Sigma_{R'}$.
	Then consider
	\begin{equation}\label{bndry_tangential_seq}
		\begin{cases}
			\Delta_{g_m} u_m^i = 0 & \text{ in } \Omega_{R'},          \\
			u_m^i = V_m^i          & \text{ on } \partial \Omega_{R'},
		\end{cases}
	\end{equation}
	for $1 \le i \le n-1$. To construct $u_m^n$, we consider
	\begin{equation}\label{bndry_normal_seq}
		\begin{cases}
			\Delta_{g_m} u_m^n = 0 & \text{ in } \Omega_{R'},          \\
			u_m^n = u^n            & \text{ on } \partial \Omega_{R'}.
		\end{cases}
	\end{equation}

	In this case the boundary data is $C^0(\partial \Omega_{R'})$. We construct solutions to \eqref{bndry_tangential_seq}-\eqref{bndry_normal_seq} using Lemma~\ref{thm:strong_solutions_c0}. By the Alexandrov maximum principle \cite[Theorem 9.1]{MR1814364}, we obtain $u_m \to u$ in $C^{0}(\Omega_{R'})$. Then, Lemma \ref{thm:local_boundary_w2p_nonhom} implies that, setting $V = u(\overline {B_r^+})$,
	\begin{equation}
		\begin{split}
			\|u^i_m - u^i\|_{W^{2,p}(V)} & \le C \Big( \| (\Delta_{g_m} - \Delta_g)u^i \|_{L^p(\Omega_{R'})}                               \\
			                             & \qquad + \|u^i_m - u^i\|_{W^{2-1/p,p}(\Sigma_{R'})} + \|u^i_m - u^i\|_{L^p(\Omega_{R'})} \Big).
		\end{split}
	\end{equation}

	The first term on the right-hand side goes to zero because $g_m \to g$ strongly in $W^{1,p}$. The second term goes to zero because the tangential harmonic coordinates converge strongly on the boundary. The third term tends to $0$ because $u_m^i \to u^i$ in $C^0(\Omega_{R'})$.

	Thus, $u_m \to u$ strongly in $W^{2,p}(B_r^+)$. In particular, for large $m$, $u_m$ restricted to $B_r^+$ is a boundary harmonic coordinate chart satisfying the $Q_m$-bounds in Definition \ref{def:boundary_harmonic_radius} with $Q_m \to Q$. Thus, $r_{bh}(z, g_m) \geq r$ for all sufficiently large $m$, establishing lower semicontinuity.

	For any fixed metric $g$, the function $z \mapsto r_{bh}(z, g)$ is $1$-Lipschitz with respect to the distance function $\operatorname{dist}_g$. Because $g_m \to g$ strongly in $W^{1,p}$, the metrics, and therefore the distance functions, converge uniformly. As a result, the sequence $z \mapsto r_{bh}(z, g_m)$ is uniformly equicontinuous on the compact set $\partial M$. By the Arzel\`a-Ascoli theorem, there is a subsequence converging uniformly. The infimum of a uniformly convergent sequence of functions converges to the infimum of the limit, namely
	\begin{equation}
		\lim_{m \to \infty} \inf_{z \in \partial M} r_{bh}(z, g_m) = \inf_{z \in \partial M} r_{bh}(z, g).
	\end{equation}
	This proves $\lim_{m \to \infty} r_{bh}(M, g_m) = r_{bh}(M, g)$, completing the proof.
\end{proof}

Uniform bounds on the interior and boundary harmonic radius imply geometric convergence. The case $\partial M = \emptyset$ is well-known; for expositions, see \cite{MR1001487}, \cite[Theorem 72]{MR2243772}, and especially \cite{MR1620864} for the case of $W^{1,p}$ metric tensor bounds. The proof for the case of manifolds with boundary is the same, so long as we assume there is an appropriate covering of a neighborhood of $\partial M$ of a definite size, which is provided by the lower bound on the global boundary harmonic radius. We state the result in our notation and show the needed modifications.

\begin{proposition}\label{thm:abstract_convergence}
	Let $(M_i, g_i)$ be a sequence of Riemannian $n$-manifolds with boundary. Fix $p > n$, $Q > 1$, and assume there exist positive constants $D_0$, $r_0$, and $c_0$ such that for all $i$, the sequence satisfies the uniform bounds:
	\begin{gather*}
		\diam(M_i) \leq D_0, \\
		r_{bh}(M_i, g_i) \geq r_0, \\
		r_h(x, g_i) \geq c_0 \dist_{g_i}(x, \partial M_i) \quad \text{for all } x \in M_i \setminus \partial M_i,
	\end{gather*}
	where $r_{bh}$ and $r_h$ denote the global boundary and interior $W^{1,p}$ harmonic radii, respectively.

	Then there exists a smooth $n$-manifold with boundary $M$ and a metric $g \in W^{1,p}(M)$ such that, passing to a subsequence, $(M_i, g_i)$ converges to $(M, g)$ in the weak $W^{1,p}$ topology and $(\partial M_i, h_i)$ converges to $(\partial M, h)$ in the weak $W^{1,p}$ topology. This convergence is realized by diffeomorphisms $F_i \colon M \to M_i$ that are uniformly bounded in $W^{2,p}(M, M_i)$, and the harmonic coordinate charts on $M_i$ subconverge weakly in $W^{2,p}$ to harmonic coordinates on $M$.
\end{proposition}

\begin{proof}
	We establish the result as a consequence of \cite[Theorem 72]{MR2243772}. We modify the conditions defining Petersen's $\|(M,g)\|_{C^{m,\alpha},r}$ norm to accommodate boundary coordinate charts, and we replace the $C^{m,\alpha}$ bounds with the $W^{k,p}$ bounds consistent with our definitions of the interior and boundary harmonic radii. Specifically, we require that each $M_i$ can be covered by a collection of charts $\varphi_s \colon \Omega_r \to U_s \subset M$ satisfying the following conditions:

	\begin{itemize}
		\item[(n1)] The domain $\Omega_r$ is either a Euclidean ball $B_r \subset \mathbb{R}^n$ (for interior charts) or an upper half-ball $B_r^+ \subset \mathbb{R}^n_+$ (for boundary charts). For boundary charts, $\varphi_s$ maps the flat portion $\Sigma_r = B_r^+ \cap \{x^n = 0\}$ to the manifold boundary $\partial M$. Furthermore, every metric ball $B(x, \frac{1}{10}e^{-Q}r)$ in $M$ is contained in some chart $U_s$.
		\item[(n2)] Petersen's bi-Lipschitz bounds on the coordinate charts are equivalent to uniform ellipticity bounds on the metric coefficients. Thus, we require $Q^{-1}\delta_{ij} \leq g_{ij} \leq Q\delta_{ij}$ on $\Omega_r$. For boundary charts, the induced boundary metric $h$ additionally satisfies $Q^{-1}\delta_{\alpha\beta} \leq h_{\alpha\beta} \leq Q\delta_{\alpha\beta}$ on $\Sigma_r$.
		\item[(n3)] The metric components satisfy the scaled $W^{k,p}$ bounds $r^{|\sigma|-n/p} \|\partial^\sigma g_{ij}\|_{L^p(\Omega_r)} \leq Q$ for all multi-indices $\sigma$ with $1 \leq |\sigma| \leq k$. For boundary charts, the induced boundary metric satisfies the corresponding fractional bounds $r^{|\sigma|-(n-1)/p} \|\partial^\sigma h_{\alpha\beta}\|_{L^p(\Sigma_r)} \leq Q$.
		\item[(n4)] The transition maps $\varphi_s^{-1} \circ \varphi_t$ satisfy a uniform $W^{k+1,p}$ bound depending only on $Q$ and $r$. For boundary charts, this bound applies to both the bulk transition maps between upper half-balls and their restrictions to the flat boundary.
	\end{itemize}
	We note our hypotheses provide a uniform scale at which the manifolds can be covered by such charts. Let $r_1 = \min(r_0/2, c_0 r_0 / 2)$. Then every point $x \in M_i$ is contained in a harmonic chart (either interior or boundary) that extends to a distance of at least $r_1$ around $x$ and satisfies the bounds from Definition \ref{def:boundary_harmonic_radius}. This implies that the sequence $M_i$ satisfies the modified (n1)-(n4) uniformly, for appropriate constants. In particular, the hypotheses of Theorem 72 are satisfied by taking $\alpha = 1 - n/p$ and $r = r_1$. Thus, using the proof of Theorem 72, inclusive of the boundary charts which satisfy the same H\"older conditions, we obtain a subsequential limit $(M,g)$, $g \in C^{0,\alpha}(M)$, with $M_i \to M$ in the $C^{0,\alpha'}$ topology, $0< \alpha' < \alpha$. The diffeomorphisms $F_i \colon M \to M_i$ are of regularity class $C^{1,\alpha}(M, M_i)$.

	With this result, and noting that the modified (n1)-(n4) conditions above are satisfied, it follows that the metric convergence is also in the weak $W^{1,p}$ topology and the diffeomorphisms are bounded in $W^{2,p}(M, M_i)$. We simply replace H\"older bounds with the corresponding Sobolev bounds at every step in the proof.

\end{proof}

\begin{remark} \label{rmk:abstract_convergence_generalizations}
	In case $\diam(M_i) \to \infty$ and the harmonic radius is bounded below on each set of fixed distance from a given base point, Proposition \ref{thm:abstract_convergence} can be modified to conclude pointed convergence. If the harmonic radius bounded is below on every compact subset in the interior of $M_i$, we can conclude local convergence on $M \backslash \partial M_i$. The modifications to the proof are straightforward.
\end{remark}

\section{Comparison Theory} \label{sec:comparison}
The purpose of this section is to develop enough comparison geometry to show that the classes $\mathcal M$, $\mathcal M_c$, and $\mathcal M_+$ are included in $\mathcal A^{bdd}$. We recall the definitions of these classes.

\begin{enumerate}[label=(\roman*)]
	\item The class $\mathcal M = \mathcal M(n, K_0, H_0, D_0, A_0, i_0)$ consists of compact, connected, Riemannian $n$-manifolds with boundary satisfying
	      \begin{gather*}
		      |\sec(M)| \leq K_0, \quad |\sec(\partial M)| \leq K_0, \quad |H| \leq H_0,\\
		      \diam(M) \leq D_0, \quad \vol_{n-1}(\partial M) \geq A_0, \quad i_b \geq i_0.
	      \end{gather*}
	\item The class $\mathcal M_+ = \mathcal M_+(n, K_0, H_0, D_0, A_0)$ consists of compact, connected, Riemannian $n$-manifolds with boundary satisfying
	      \begin{gather*}
		      |\sec(M)| \leq K_0, \quad |\sec(\partial M)| \leq K_0, \quad 1/H_0 < |H| \leq H_0,\\
		      \diam(M) \leq D_0, \quad \vol_{n-1}(\partial M) \geq A_0.
	      \end{gather*}
	\item The class $\mathcal{M}_c = \mathcal{M}_c(n, K_0, S_0, D_0, v_0)$ consists of smooth, compact, connected Riemannian $n$-manifolds with boundary satisfying
	      \begin{gather*}
		      |\sec(M)| \leq K_0, \quad 0 \leq S \leq S_0, \\
		      \diam(M) \leq D_0, \quad \vol(M) \geq v_0.
	      \end{gather*}
	\item The class $\mathcal A^{\mathrm{bdd}} = \mathcal{A}^{\mathrm{bdd}}(n, K_0, i_0, D_0, H_0)$ consists of smooth, compact, connected Riemannian $n$-manifolds with boundary satisfying
	      \begin{gather*}
		      \|\Ric_M\|_{L^{\infty}(M)} \leq K_0, \quad \|\Ric_{\partial M}\|_{L^{\infty}(\partial M)} \leq K_0, \\
		      \diam(M) \leq D_0, \quad |H| \leq H_0, \\
		      i_M \geq i_0, \quad i_{\partial M} \geq i_0, \quad i_b \geq i_0.
	      \end{gather*}
\end{enumerate}

Let us suppose throughout that $(M,g)$ is a smooth, connected, compact Riemannian manifold with nonempty boundary. Write $\tau(x) = \dist(x, \partial M)$ and note that $\tau$ is smooth off of the cut locus of $\partial M$. Define
\begin{equation*}
	\nu = \grad \tau
\end{equation*}
so that $\nu$ is the inward-pointing unit normal to $\partial M$. Let $S(X) = -\nabla_X \nu$ denote the shape operator, where $\nabla$ is the usual covariant derivative, and write $g(S(X), Y)$ for the second fundamental form. The mean curvature $H$ on $\partial M$ is $H := \frac{1}{n-1}\tr S$. We begin with the basic observation that the conditions defining $\mathcal M_+$ and $\mathcal M$ imply uniform bounds on the second fundamental form. Recall the twice-traced Gauss equation for the boundary $\partial M \subset M$,
\begin{equation} \label{eq:twice_traced_gauss}
	R_{\partial M} = R_M - 2 \Ric_M(\nu, \nu) + (n-1)^2 H^2 - |S|^2,
\end{equation}
where $R_{\partial M}$ and $R_M$ are the scalar curvatures of $\partial M$ and $M$, respectively. Rearranging this equation, we obtain
\begin{equation*}
	|S|^2 = (n-1)^2 H^2 + R_M - 2 \Ric_M(\nu, \nu) - R_{\partial M}.
\end{equation*}
For any element $(M, g)$ in $\mathcal M$ or $\mathcal M_+$, the sectional curvatures of $M$ and $\partial M$ are bounded in absolute value by $K_0$, which implies uniform bounds on the scalar curvatures and the Ricci curvature. Since we also have $|H| \le H_0$, we obtain a uniform bound on the second fundamental form $|S|$ of the boundary. Similarly, if $(M,g) \in \mathcal M_c$, the Gauss equation implies a uniform bound on the sectional curvature of the boundary $\partial M$.

To make further observations, let us first write $\foc(\partial M)$ for the focal locus distance of $\partial M$. We note the following classical bounds (see, for instance, \cite{warner1966extensions}).
\begin{lemma}[Warner's focal locus bound]\label{prop:focal_bound}
	Suppose the principal curvatures of $\partial M$ are bounded above by $\lambda$, and the sectional curvature of $M$ is bounded above by $K$. Let $t_0$ be the smallest positive solution to
	\begin{align*}
		\cot(\sqrt{K} t)   & = \frac{\lambda}{\sqrt{K}}  &  & \textrm{if } K > 0,  \\
		t                  & = \frac{1}{\lambda}         &  & \textrm{if } K  = 0, \\
		\coth(\sqrt{-K} t) & = \frac{\lambda}{\sqrt{-K}} &  & \textrm{if } K < 0.
	\end{align*}
	If no such positive solution exists, set $t_0 = \infty$. Then $\foc(\partial M) \geq t_0$.

	Suppose the principal curvatures of $\partial M$ are bounded below by $\Lambda$, and the sectional curvature of $M$ is bounded below by $\kappa$. Let $t_1$ be the smallest positive solution to
	\begin{align*}
		\cot(\sqrt{\kappa} t)   & = \frac{\Lambda}{\sqrt{\kappa}}  &  & \textrm{if } \kappa > 0,                                         \\
		t                       & = \frac{1}{\Lambda}              &  & \textrm{if } \kappa  = 0,                                        \\
		\coth(\sqrt{-\kappa} t) & = \frac{\Lambda}{\sqrt{-\kappa}} &  & \textrm{if } \kappa < 0 \textrm{ and } \Lambda > \sqrt{-\kappa}.
	\end{align*}
	If a finite positive solution $t_1$ exists, then $\foc(\partial M) \leq t_1$.
\end{lemma}

Write $\nu(\partial M)$ for the normal bundle of $\partial M$ and define the normal exponential map $\exp_{\nu}:\nu(\partial M) \to M$. The boundary injectivity radius $i_b(M)$ is defined as the supremum of $t$ such that $\exp_{\nu}$ is a diffeomorphism on $\partial M \times[0,t)$. It is well-known that
\begin{align*}
	i_b(M) \geq \min\{\foc(\partial M), \cut(\partial M)\},
\end{align*}
where $\cut(\partial M)$ is the cut locus distance to $\partial M$, defined as the infimum of the distance to the first cut point along normal geodesics starting from $\partial M$. The following result is also classical; we include a proof here to show dependency on constants (see also  \cite[Lemma 2.4]{BoundaryMetric3mfld}).

\begin{lemma}\label{lem:cutlocuslemma}
	Suppose $(M,g) \in \mathcal M_{+}$. Then
	\begin{equation*}
		i_b(M) \geq \min\left\{t_0, \frac{2}{H_0K_0}\right\},
	\end{equation*}
	where $t_0$ is the focal bound defined in Lemma \ref{prop:focal_bound}. Thus $\mathcal M_+ \subset \mathcal M$. If, in addition, $\Ric(g) \geq 0$, then $i_b(M) \geq t_0$.
\end{lemma}

\begin{proof}
	By Proposition \ref{prop:focal_bound}, $\foc(\partial M) \geq t_0$. Thus it suffices to show that if $\cut(\partial M) < t_0$, then
	\begin{align*}
		\cut(\partial M) \geq \frac{2}{H_0K_0}.
	\end{align*}
	Assume $\cut(\partial M) < t_0$. We may choose an arclength-parametrized geodesic $\gamma: [0,l] \to M$ realizing the cut locus distance $l = \cut(\partial M)$. Thus, $\gamma$ is a minimizing geodesic from $p \in \partial M$ to $q \in \partial M$ that is orthogonal to $\partial M$ at both endpoints, and $\im(\gamma) \cap \partial M = \{p,q\}$.

	Consider the index form for vector fields $V, W$ along $\gamma$ that are orthogonal to $\gamma'$:
	\begin{align*}
		I(V,W) = \int_0^l \left( g(\nabla_{\gamma'}V,\nabla_{\gamma'}W) - g(R(\gamma',V)W,\gamma') \right) dt - g(S(V),W)\big|_0^l.
	\end{align*}
	Because $l = \cut(\partial M) < t_0 \leq \foc(\partial M)$, there are no focal points along $\gamma$. Consequently, the index form is positive semi-definite, meaning $I(V,V) \geq 0$ for any such vector field $V$.

	Choose an orthonormal basis $\{e_i\}_{i=1}^{n-1}$ for $T_p \partial M$ and let $V_i$ be the parallel translation of $e_i$ along $\gamma$. Taking the trace of the index form over this basis gives
	\begin{align*}
		0 \leq \sum_{i=1}^{n-1} I(V_i,V_i) & = -\int_0^l \Ric(\gamma',\gamma')\,dt -(n-1) (H(p) + H(q)) \\
		                                   & \leq (n-1)lK_0 - \frac{2(n-1)}{H_0}.
	\end{align*}
	Rearranging this inequality gives $l \geq \frac{2}{H_0K_0}$, establishing the first assertion.

	If we additionally assume $\Ric(g) \geq 0$, the trace inequality becomes:
	\begin{align*}
		0 \leq -\int_0^l \Ric(\gamma',\gamma')\,dt - (H(p) + H(q)) \leq -\frac{2}{H_0} < 0,
	\end{align*}
	which is a contradiction. Therefore, if $\Ric(g) \geq 0$, it is impossible for $\cut(\partial M)$ to be strictly less than $t_0$, implying $i_b(M) \geq t_0$.
\end{proof}
Next, we control the mean curvature on level sets near the boundary, using an argument parallel to the asymptotic case discussed in \cite[Lemma 3.2.2]{MR2096795}.
\begin{lemma}\label{lem:optimal_mean_curvature}
	Suppose $(M, g)$ is a Riemannian $n$-manifold with boundary satisfying the Ricci curvature lower bound $\Ric \ge -(n-1)K_0$ for some $K_0 \ge 0$. Let $i_b > 0$ denote the boundary injectivity radius. Let $\Sigma_t$ denote the level set of points at distance $t$ from $\partial M$, and let $H(t)$ denote the mean curvature of $\Sigma_t$. Then for all $0 < t < i_b/2$, we have
	\[
		H(t) \le \frac{1}{c_1(i_b - t)},
	\]
	where $c_1 = 1$ if $K_0 = 0$, and if $K_0 > 0$,
	\[
		c_1 = \frac{\sqrt{1 + 4K_0 i_b^2} - 1}{2K_0 i_b^2}.
	\]
\end{lemma}

\begin{proof}
	Write $\nu$ for the inward normal to $\Sigma_t$, let $S(t) = -\nabla \nu$ be the shape operator, and define $h(t) = -\operatorname{tr}(S) = -(n-1)H(t)$. By the standard Riccati equation along normal geodesics, $h'(t) + \operatorname{tr}(S^2) = -\Ric(\nu, \nu)$. Using the Cauchy-Schwarz inequality $\operatorname{tr}(S^2) \ge \frac{h(t)^2}{n-1}$ and the curvature bound $\Ric(\nu, \nu) \ge -(n-1)K_0$, we obtain the differential inequality:
	\begin{equation} \label{eq:riccati_ineq}
		h'(t) \le -\frac{h(t)^2}{n-1} + (n-1)K_0.
	\end{equation}

	If $K_0 = 0$, \eqref{eq:riccati_ineq} reduces to $h'(t) \le -\frac{h(t)^2}{n-1}$. Integrating this from $t$ to $s$ (for $t < s < i_b$) yields $\frac{1}{h(s)} - \frac{1}{h(t)} \ge \frac{s - t}{n-1}$. Since boundary injectivity ensures no focal points occur before $i_b$, we have $\limsup_{s \to i_b} \frac{1}{h(s)} \le 0$. Taking the limit as $s \to i_b$ gives $h(t) \ge -\frac{n-1}{i_b - t}$, which corresponds to $c_1 = 1$.

	Assume $K_0 > 0$. By definition, $c_1 \in (0, 1)$ is the unique positive root of the equation $\frac{c_1^2}{1-c_1} = \frac{1}{K_0 i_b^2}$. For any $t \in (0, i_b/2)$, we have $(i_b - t)^2 < i_b^2$, which implies the strict inequality:
	\begin{equation} \label{eq:strict_c1}
		\frac{1}{K_0(i_b - t)^2} > \frac{c_1^2}{1-c_1}.
	\end{equation}

	We want to show $h(t) \ge -\frac{n-1}{c_1(i_b - t)}$. Assume for contradiction there exists $t_0 \in (0, i_b/2)$ such that $h(t_0) < -\frac{n-1}{c_1(i_b - t_0)} < 0$. Squaring this yields $h(t_0)^2 > \frac{(n-1)^2}{c_1^2(i_b - t_0)^2}$. Rearranging \eqref{eq:strict_c1} at $t=t_0$ and substituting this bound gives:
	\[
		(n-1)K_0 < \frac{1-c_1}{c_1^2(i_b - t_0)^2}(n-1) < \frac{1-c_1}{n-1}h(t_0)^2.
	\]
	Plugging this into \eqref{eq:riccati_ineq} at $t=t_0$, we find
	\[
		h'(t_0) < -\frac{h(t_0)^2}{n-1} + \frac{1-c_1}{n-1}h(t_0)^2 = -\frac{c_1}{n-1}h(t_0)^2.
	\]
	Let $B(t) = -\frac{n-1}{c_1(i_b - t)}$. Notice that $B'(t) = -\frac{n-1}{c_1(i_b - t)^2}$. As long as $h(t) \le B(t)$, the same logic guarantees $h'(t) < -\frac{c_1}{n-1}h(t)^2 \le -\frac{c_1}{n-1}B(t)^2 = B'(t)$. Since $h(t_0) < B(t_0)$, it follows that $h(t)$ decreases strictly faster than $B(t)$, preserving the strict inequality $h'(t) < -\frac{c_1}{n-1}h(t)^2$ for all $t \in[t_0, i_b)$.

	Dividing by $h(t)^2$ and integrating from $t_0$ to $t$ we obtain
	\[
		-\frac{1}{h(t)} < -\frac{1}{h(t_0)} - \frac{c_1}{n-1}(t - t_0).
	\]
	Because $h(t_0) < B(t_0)$, we can write $-\frac{1}{h(t_0)} = \frac{c_1(i_b - t_0)}{n-1} - \epsilon$ for some $\epsilon > 0$. Substituting this into the integrated inequality gives:
	\[
		-\frac{1}{h(t)} < \frac{c_1(i_b - t)}{n-1} - \epsilon.
	\]
	As $t \to i_b$, the right-hand side approaches $-\epsilon < 0$. However, since $h(t)$ is strictly negative, $-\frac{1}{h(t)}$ must remain strictly positive as long as $h(t)$ is finite. The inequality forces $-\frac{1}{h(t)}$ to hit zero at some $t^* < i_b$, meaning $h(t) \to -\infty$. This implies a focal point forms strictly before distance $i_b$, contradicting the definition of the boundary injectivity radius. Thus, our assumption is false, and we conclude $h(t) \ge -\frac{n-1}{c_1(i_b - t)}$, meaning $H(t) = -\frac{h(t)}{n-1} \le \frac{1}{c_1(i_b - t)}$ for all $0 < t < i_b/2$.
\end{proof}

For any $M \in \mathcal M$, lemma \ref{lem:optimal_mean_curvature} provides a uniform estimate on the mean curvature of level sets near the boundary. As an application, we bound the volume of small cylinders in $M$ with base $\Sigma_0 \subset \partial M$.
\begin{lemma}\label{lem:isoperimetric}
	Suppose $(M,g) \in \mathcal{M}$. Fix $z \in \partial M$, $s < \diam(\partial M)$, $0 \le t_1 \le t_2 < i_b/2$, and write $\Sigma_0 := B_{z}^{\partial M}(s)$, the ball of radius $s$ within $\partial M$. Define the cylinder
	\begin{equation*}
		\operatorname{Cyl}(\Sigma_0,t_1,t_2) := \{ \exp_{\nu}(p,t) \textrm{ : } t_1 \leq t \leq t_2\textrm{, } p \in \Sigma_0\}.
	\end{equation*}
	Then there exists a constant $a_0 > 0$, depending only upon the constants defining $\mathcal{M}$, so that
	\begin{equation*}
		\vol(\operatorname{Cyl}(\Sigma_0,t_1,t_2)) \geq a_0s^{n-1}(t_2 - t_1).
	\end{equation*}
\end{lemma}

\begin{proof}
	First, a standard volume comparison on the boundary manifold $\partial M$ (which satisfies $\sec(\partial M) \ge -K_0$) implies that
	\begin{equation} \label{eq:boundary_vol_comp}
		\vol_{n-1}(\Sigma_0) \geq c s^{n-1}
	\end{equation}
	for some constant $c > 0$ depending only on $n$, $K_0$, $\diam(\partial M)$, and $A_0$. We note here that, because the extrinsic diameter of $M$ is bounded by $D_0$ and the second fundamental form is uniformly bounded in terms of the constants defining $\mathcal M$, a result of Wong \cite[Thm 1.1]{MR2379775} guarantees a bound on $\diam(\partial M)$ that similarly depends only on the constants defining $\mathcal M$.

	Let $\Sigma_{r} = \{ \exp_{\nu}(p,r) \textrm{ : } p \in \Sigma_0 \}$ be the level set at distance $r$. Since $\vol(\operatorname{Cyl}(\Sigma_0,t_1,t_2)) = \int_{t_1}^{t_2}\!\vol_{n-1}(\Sigma_r)\,dr$, it suffices to uniformly bound $\vol_{n-1}(\Sigma_r)$ from below.

	By Lemma \ref{lem:optimal_mean_curvature}, the mean curvature of the level sets satisfies $H(r) \leq \frac{1}{c_1(i_b -r)}$ for all $0 < r < i_b/2$. Write $A(r) := \vol_{n-1}(\Sigma_r)$ and $A_0 := \vol_{n-1}(\Sigma_0)$. The first variation of area is given by
	\begin{equation*}
		A'(r) = -(n-1)\int_{\Sigma_r}\! H d\mu_r,
	\end{equation*}
	where $d\mu_r$ is the volume form on $\Sigma_r$. Applying the mean curvature bound, we obtain a differential inequality
	\begin{equation*}
		A'(r) \geq -\int_{\Sigma_r}\frac{n-1}{c_1(i_b -r)} d\mu_r = -\frac{n-1}{c_1(i_b -r)}A(r).
	\end{equation*}
	Dividing by $A(r)$ and integrating from $0$ to $r$ gives
	\begin{equation*}
		\ln\left(\frac{A(r)}{A_0}\right) \geq \frac{n-1}{c_1} \ln\left(\frac{i_b - r}{i_b}\right) \implies A(r) \geq A_0 \left( \frac{i_b -r}{i_b} \right)^{\frac{n-1}{c_1}}.
	\end{equation*}
	Because we restrict to $r \le t_2 < i_b/2$, we have $\frac{i_b - r}{i_b} > \frac{1}{2}$. Thus, $A(r) \ge 2^{-\frac{n-1}{c_1}}A_0$.

	Integrating this lower bound from $t_1$ to $t_2$ and applying \eqref{eq:boundary_vol_comp}, we obtain
	\begin{equation*}
		\vol(\operatorname{Cyl}(\Sigma_0,t_1,t_2)) = \int_{t_1}^{t_2} A(r) \, dr \geq A_0 \left(\frac{1}{2}\right)^{\frac{n-1}{c_1}} (t_2 - t_1) \geq c \left(\frac{1}{2}\right)^{\frac{n-1}{c_1}} s^{n-1} (t_2 - t_1),
	\end{equation*}
	completing the proof.
\end{proof}
We use this result to show that elements of $\mathcal M$ cannot undergo interior volume collapse.
\begin{proposition}\label{prop:interior_volume}
	There exists $v_0 > 0$, depending only on the constants that determine $\mathcal M$, so that for every $(M,g) \in \mathcal M$, $x \in M \setminus \partial M$, and $r < \dist(x,\partial M)$, there holds
	\begin{equation*}
		\vol(B_r(x)) \geq v_0r^n.
	\end{equation*}
\end{proposition}

\begin{proof}
	Let $r_x = \dist(x,\partial M)$. Choose a unit-speed, distance-minimizing geodesic $\gamma$ from $\partial M$ to $x$, such that $\gamma(0) = z \in \partial M$ and $\gamma(r_x) = x$, meeting $\partial M$ orthogonally at $z$. Write $r_y = \min(r_x, i_b/2)$ and define $y = \gamma(r_y)$. By standard Jacobi field estimates (since the sectional curvatures of $M$ and the principal curvatures of $\partial M$ are bounded), the differential of the normal exponential map is uniformly bounded on $t \in [0, i_b/2]$. Thus, there exists a constant $C > 0$ such that for any $t \leq i_b/2$, if $\dist_{\partial M}(z', z) \leq \epsilon$, the distance in $M$ satisfies $\dist(\exp_{\nu}(z',t), \exp_{\nu}(z, t)) \leq C\epsilon$.

	Define the boundary ball $\Sigma_0 := B^{\partial M}_{r_y / 4C}(z)$. By the triangle inequality and the distance distortion bound, the cylinder over $\Sigma_0$ satisfies the inclusion
	\begin{equation} \label{eq:cylinder_inclusion}
		\operatorname{Cyl}(\Sigma_0, 3r_y/4, r_y) \subset B_{r_y/2}(y) \subset B_{r_x}(x).
	\end{equation}
	To verify this inclusion, note that for any point $w = \exp_\nu(z',t)$ in the cylinder, its distance to $\gamma(t) = \exp_\nu(z,t)$ in $M$ is at most $r_y/4$, and $t$ is within $r_y/4$ of $r_y$. Thus, its distance to $y = \gamma(r_y)$ is at most $r_y/2$. Since $y$ lies on the minimizing geodesic from $\partial M$ to $x$ at distance $r_x - r_y$ from $x$, the ball $B_{r_y/2}(y)$ is strictly contained in $B_{r_x}(x)$.

	Applying Lemma \ref{lem:isoperimetric} to this cylinder with $s = r_y / 4C$, $t_1 = 3r_y/4$, and $t_2 = r_y$, we obtain:
	\begin{equation*}
		\vol(B_{r_x}(x)) \geq \vol(\operatorname{Cyl}(\Sigma_0, 3r_y/4, r_y)) \geq a_0 \left(\frac{r_y}{4C}\right)^{n-1} \left(\frac{r_y}{4}\right) = a_1 r_y^n,
	\end{equation*}
	where $a_1 = a_0 / (4^n C^{n-1})$. Since $r_y = \min(r_x, i_b/2)$, we have $\vol(B_{r_x}(x)) \geq a_1 \min(r_x^n, (i_b/2)^n)$.

	Finally, for any $r < r_x$, we apply the Bishop-Gromov volume comparison theorem to the balls $B_r(x)$ and $B_{r_x}(x)$. Because the Ricci curvature is bounded below by $-(n-1)K$, this implies:
	\begin{equation*}
		\frac{\vol(B_r(x))}{\vol_K(r)} \geq \frac{\vol(B_{r_x}(x))}{\vol_K(r_x)},
	\end{equation*}
	where $\vol_K$ denotes the volume of a ball in the space form of constant curvature $-K$. Substituting our lower bound for $\vol(B_{r_x}(x))$, we obtain
	\begin{equation*}
		\vol(B_r(x)) \geq a_1 \min(r_x^n, (i_b/2)^n) \frac{\vol_K(r)}{\vol_K(r_x)}.
	\end{equation*}
	Because $r_x \leq \diam(M) \leq D$, the ratio $\min(r_x^n, (i_b/2)^n)/\vol_K(r_x)$ is uniformly bounded from below on $(0, D]$. Indeed, as $r_x \to 0$, this ratio continuously approaches $1/\omega_n$ (where $\omega_n$ is the volume of the Euclidean unit ball), and for $r_x > 0$ it remains strictly positive, thus attaining a minimum that depends only on $i_b$, $K_0$, $n$, and $D_0$. Furthermore, the space form volume satisfies $\vol_K(r) \geq \omega_n r^n$. Combining these facts yields a uniform constant $v_0 > 0$, depending only on $a_1$, $i_b$, $K_0$, $n$, and the diameter bound $D_0$, such that $\vol(B_r(x)) \geq v_0 r^n$.
\end{proof}

Finally, we tie together the results of this section to identify subclasses of $\mathcal A^{bdd}$.

\begin{proposition} \label{cor:strata_inclusion}
	It holds that $\mathcal M_+ \subset \mathcal M \subset \mathcal A^{bdd}$.
\end{proposition}

\begin{proof}
	We've already shown in Lemma \ref{lem:cutlocuslemma} that $\mathcal M_+ \subset \mathcal M$. Considering $\mathcal{M}$, the diameter, mean curvature, and boundary injectivity radius bounds are assumed by definition, and the sectional curvature bounds imply corresponding Ricci curvature bounds. A lower bound on the injectivity radius of the intrinsic boundary follows from standard Cheeger-Gromov theory.

	The more substantive requirement is to establish a uniform lower bound on the interior injectivity radius. Thus, suppose $(M,g) \in \mathcal M$ and $x \in M \setminus \partial M$. By Proposition \ref{prop:interior_volume}, there exists a uniform constant $v_0 > 0$ such that $\vol(B_r(x)) \ge v_0 r^n$ for all $r < d := \dist(x, \partial M)$. Consider the metric $\tilde g = d^{-2} g$, in which the distance to the boundary is normalized to $1$. The unit ball $B_1(x, \tilde g) $ satisfies the rescaled sectional curvature bound $|\sec(\tilde g)| \le K_0 d^2 \le K_0 D_0^2$ (where $D_0$ is the uniform diameter bound) and the volume bound $\operatorname{vol}_{\tilde g}(B_1(x, \tilde g)) \ge v_0$. By the local injectivity radius estimate of Cheeger, Gromov, and Taylor \cite[Theorem 4.7]{MR658471}, these uniform bounds on the curvature and volume guarantee a uniform lower bound on the injectivity radius at $x$ in the rescaled metric,
	\begin{equation*}
		\operatorname{inj}_{\tilde g}(x) \ge c_0(n, K_0 D_0^2, v_0) > 0.
	\end{equation*}
	With the original metric $g$, we obtain a scale-invariant lower bound on the interior injectivity radius, namely
	\begin{equation*}
		\operatorname{inj}_g(x) \ge c_0 \operatorname{dist}_g(x, \partial M).
	\end{equation*}
\end{proof}

\begin{proposition}\label{thm:kodani_improvement}
	It holds that $\mathcal{M}_{c} \subset \mathcal{A}^{\mathrm{bdd}}$.
\end{proposition}

\begin{proof}
	To establish inclusion in $\mathcal{A}^{\mathrm{bdd}}$, we must verify the uniform bounds on the Ricci curvatures, the mean curvature, and the three injectivity radii $i_M$, $i_{\partial M}$, and $i_b$. As noted at the beginning of Section \ref{sec:comparison}, the Gauss equation and the bounds on the sectional curvature and second fundamental form immediately imply uniform $L^\infty$ bounds on the Ricci curvature of both the interior and the boundary and the mean curvature.

	Because the extrinsic diameter is bounded by $D_0$ and the second fundamental form is uniformly bounded by $S_0$, a result of Wong \cite[Thm 1.1]{MR2379775} guarantees a uniform upper bound on the intrinsic diameter of the boundary, $\operatorname{diam}(\partial M) \le D_{\partial}$.

	Uniform lower bounds for both the interior injectivity radius $i_M \ge i_0$ and the boundary injectivity radius $i_b \ge i_0$ are established by Kodani in \cite[Proposition 6.1]{MR1065204}.

	It remains only to bound the intrinsic boundary injectivity radius $i_{\partial M}$. The Heintze-Karcher theorem also provides a uniform constant $C = C(n, K_0, S_0, D_0)$ such that the volume of $M$ is controlled by the $\vol_{n-1}$ of its boundary: $\operatorname{vol}(M) \le C \vol_{n-1}(\partial M)$. The interior volume non-collapsing $\operatorname{vol}(M) \ge v_0$ therefore forces a uniform lower bound on the boundary $\vol_{n-1}$: $\vol_{n-1}(\partial M) \ge v_0 / C > 0$. Thus, standard Cheeger-Gromov theory guarantees a uniform lower bound on its intrinsic injectivity radius, $i_{\partial M} \ge i_0$.
\end{proof}

\section{Compactness Theorems}\label{section_compactness}
Recall that $\mathcal A^{bdd}$ is the class of smooth Riemannian manifolds with boundary $(M, g)$ satisfying
\begin{gather*}
	\|\Ric_M\|_{L^{\infty}(M)} \leq K_0, \quad \|\Ric_{\partial M}\|_{L^{\infty}(\partial M)} \leq K_0, \\
	\|H\|_{L^{\infty}(\partial M)} \leq H_0, \quad \diam(M) \leq D_0,\\
	i_M \geq i_0, \quad i_{\partial M} \geq i_0, \quad i_b \geq i_0.
\end{gather*}
The difference between this class and the class $\mathcal A$ studied in \cite{MR2096795} is that the Lipschitz condition on $H$ has been relaxed to an $L^{\infty}$ bound. This section is devoted to showing $\mathcal A^{bdd}$ is precompact in the $W^{1,p}$ topology, any $p > n$. The main argument is Theorem \ref{precompact}, while much of the analysis is found in Proposition \ref{lem:local_strong_convergence}, which uses the results of Section \ref{defs} to strengthen the abstract convergence result Proposition \ref{thm:abstract_convergence} under the presence of curvature bounds.

\begin{theorem}\label{precompact}
	The class $\mathcal{A}^{\mathrm{bdd}}(n, K_0, i_M, i_0, D_0, H_0)$ is precompact in the $W^{1,p}$ topology for any $p > n$. Consequently, it is also precompact in the $C^\alpha$ topology for any $0 < \alpha < 1$. Specifically, if $(M_k, g_k) \in \mathcal{A}^{\mathrm{bdd}}$ is a sequence of Riemannian manifolds with boundary, then there exists a smooth manifold with boundary $M$ equipped with a $W^{1,p}$ metric tensor $g$, so that after passing to a subsequence, there are $W^{2,p}$ diffeomorphisms $F_k: M \to M_k$ with $F_k^*g_k \to g$ in the $W^{1,p}$ topology on $M$.
\end{theorem}

\begin{proof}
	First, we show that $\mathcal A^{bdd}$ is weakly precompact in $W^{1,p}$, any $p > n$. By Proposition \ref{thm:abstract_convergence}, it suffices to prove there exist uniform constants $r_0 > 0$ and $c_0 > 0$ such that for all $(M, g) \in \mathcal{A}^{\mathrm{bdd}}$, the boundary harmonic radius satisfies $r_{bh}(M, g) \ge r_0$, and the interior harmonic radius satisfies $r_h(x, g) \ge c_0 \dist_{g}(x, \partial M)$ for all $x \in M \setminus \partial M$.

	We will focus on the boundary harmonic radius bound, as the interior case follows from \cite{MR1074481} without any modification (in fact, stronger, $C^{1,\alpha}$ convergence holds good away from the boundary.) See also \cite{MR2096795} for another exposition of the interior case.

	Assume for contradiction that the uniform boundary harmonic radius bound fails. Then there exists a sequence of manifolds $(M_k, \tilde{g}_k) \in \mathcal{A}^{\mathrm{bdd}}$ such that the global boundary harmonic radius satisfies $r_k := r_{bh}(M_k, \tilde{g}_k) \to 0$ as $k \to \infty$.

	We rescale the metrics by setting $g_k = r_k^{-2} \tilde{g}_k$. Under this rescaling, the boundary harmonic radius normalizes to $r_{bh}(M_k, g_k) = 1$. We then choose points $p_k \in \partial M_k$ so that $r_{bh}(p_k, g_k) = 1$. The geometric bounds scale as follows:
	\begin{gather*}
		\|\Ric_{M_k}\|_{L^\infty(g_k)} = r_k^2 \|\Ric_{M_k}\|_{L^\infty(\tilde{g}_k)} \le r_k^2 K_0 \to 0, \\
		\|\Ric_{\partial M_k}\|_{L^\infty(g_k)} = r_k^2 \|\Ric_{\partial M_k}\|_{L^\infty(\tilde{g}_k)} \le r_k^2 K_0 \to 0, \\
		\|H\|_{L^\infty(g_k)} = r_k \|H\|_{L^\infty(\tilde{g}_k)} \le r_k H_0 \to 0, \\
		i_{M_k}(g_k) = r_k^{-1} i_{M_k}(\tilde{g}_k) \ge r_k^{-1} i_0 \to \infty, \\
		i_{\partial M_k}(g_k) = r_k^{-1} i_{\partial M_k}(\tilde{g}_k) \ge r_k^{-1} i_0 \to \infty, \\
		i_b(g_k) = r_k^{-1} i_b(\tilde{g}_k) \ge r_k^{-1} (2i_0) \to \infty, \\
		\diam(g_k) = r_k^{-1} \diam(\tilde{g}_k) \le r_k^{-1} D_0 \to \infty.
	\end{gather*}
	By Remark \ref{rmk:abstract_convergence_generalizations}, the sequence of pointed manifolds $(M_k, g_k, p_k)$ subconverges in the pointed weak $W^{1,p}$ topology to a limit $(M, g, p)$. (Actually, to apply Remark \ref{rmk:abstract_convergence_generalizations}, we still need to argue that $r_h(x, g_k) \geq c_0 \dist_{g_k}(x, \partial M_k)$ for all $x \in M_k \backslash \partial M_k$. This argument is the same as the interior case already proved in \cite{MR1074481}.)

	By the definition of pointed weak $W^{1,q}$ convergence, for any compact subset $V \subset M$, the pull-back metrics (which we continue to denote by $g_k$ for simplicity) subconverge weakly in $W^{1,p}$ to $g$ on $V$. Suppose $z \in V \cap \partial M$. Choose $R > 0$ so that there exists a $g$-harmonic boundary chart $u \colon U \to B_R^+$  and $g_k$-harmonic boundary charts $u_k \colon U_k \to B_R^+$ centered at $z$, with the
	transition maps converging to the identity weakly in $W^{2,p}(B_R^+)$. Such charts are explicitly constructed in the proof of Proposition \ref{weaktostrong}. We apply Proposition \ref{lem:local_strong_convergence} below to conclude that $(g_k)_{ij} \to g_{ij}$ strongly in $W^{1,q}_{\mathrm{loc}}(B^+_{R})$. As in the proof of Proposition \ref{weaktostrong}, the $g_k$ converge strongly in $W^{1,q}_{\mathrm{loc}}$ on each interior chart (in fact, the interior convergence is stronger, in $C^{1,\alpha}_{loc}(M \setminus \partial M))$. Patching these local limits together, we conclude that $(M_k, g_k, p_k)$ converges in $W^{1,q}$, any $n < q < \infty$, to its limit $(M, g, p)$ and that $g \in W^{1,q}(M)$ for all $q > n$.

	The limiting manifold satisfies
	\begin{gather*}
		\Ric(g) = 0, \quad
		\Ric_{\partial M}(h) = 0, \quad H =0.
	\end{gather*}
	Interior $C^{\infty}$ regularity is well known, coming from the equation $\Ric(g) = 0$. To establish boundary regularity, we apply \cite[Theorem 2.1]{MR2096795} to conclude $g \in C^{1,\alpha}$. Standard elliptic theory then bootstraps the regularity to $C^{\infty}$. By Lemma \ref{lem:anderson_w1p} below, we conclude that the limit is the Euclidean upper half-space.
	Because the boundary harmonic radius is continuous under strong $W^{1,p}$ convergence (Proposition \ref{thm:harmonic_radius_continuity}), we must have $r_{bh}(p,g) = \lim_{k \to \infty} r_{bh}(p_k, g_k) = 1$, contradicting the fact that the limit is the Euclidean upper half-space.

	This establishes uniform lower bounds for the interior and boundary harmonic radii. Together with Proposition \ref{thm:abstract_convergence}, this shows that any sequence $(M_k, g_k) \in \mathcal A^{bdd}$ has a weakly convergent $W^{1,p}$ subsequence, any $p > n$. Finally, applying Proposition \ref{weaktostrong}, we extract a strongly convergent subsequence, finishing the proof.

\end{proof}

\begin{proposition}[Weak-to-Strong $W^{1,p}$ Convergence] \label{lem:local_strong_convergence}
	Suppose $M$ is a smooth manifold with boundary and $g_k$ is a sequence of $W^{1,p}$ Riemannian metrics with $g_k \rightharpoonup g$ weakly in $W^{1,p}(M)$ and intrinsic boundary metrics $h_k \to h$ weakly in $W^{1,p}(\partial M)$. Fix $z \in \partial M$ and suppose there are open sets $U_k \ni z$ and $R > 0$ with $u_k : U_k \to B_R^+$ being a boundary harmonic coordinate chart for $g_{k}$ centered at $z$, where $B_R^+ = B_R \cap \{x^n > 0\}$ is the open upper half-ball. Write $\Sigma_R = \overline{B_R^+} \cap \{x^n = 0\}$ for the tangential boundary portion. Suppose that $U \ni z$ is a relatively open set, $u: U \to B_R^+$ is a boundary harmonic coordinate chart for $g$, and $u_k \to u$ weakly in $W^{2,q}(B_R^+)$, any $q > n$. Assume further that there exist constants $K, H_0$ such that
	\[
		\|\operatorname{Ric}(g_k)\|_{L^\infty(B_R^+)} \leq K, \quad \|\operatorname{Ric}(h_k)\|_{L^\infty(\Sigma_R)} \leq K, \quad \|H_{g_k}\|_{L^\infty(\Sigma_R)} \leq H_0.
	\]
	Then for any $r < R$, $g_k \to g$ in $W^{1,q}(B_r^+)$, any $q > n$.
\end{proposition}

\begin{proof}
	Using Proposition \ref{thm:boundary_regularity} (see also Corollary \ref{prop:local_neumann_l_s}), we conclude that $\{g_k\}$ are uniformly bounded in $W^{1,q}(B_R^+)$ for any $q < \infty$. In addition, we have that the gradient of $g_k$ admits a trace on $\Sigma_R$ that is bounded in $L^q(\Sigma_R)$, any $q > n$. A first conclusion is that $g_k \to g$ in $C^{0,\alpha}(B_R^+)$, $0 <\alpha <1$.  We will establish strong convergence in $W^{1,q}(B_r^+)$ by further analyzing the elliptic system \eqref{eq:tangential_interior}--\eqref{eq:normal_boundary_2} on $B_R^+$ in local boundary harmonic coordinates $u_k$ for $g_k$ and $u$ for $g$. From our hypothesis, the optimal regularity property of harmonic coordinates, and the further bounds just established, in the coordinates $u$ on $B_R^+$ we have that $g_k \rightharpoonup g$ in $W^{1,q}(B_R^+)$. In fact, due to $u_k \rightharpoonup u$ in $W^{2,q}(B^+_R)$ we have $(g_k)_{ij} \rightharpoonup g_{ij}$ in $W^{1,q}(B_R^+)$ where $(g_k)_{ij}$ is expressed in the $u_k$ coordinates and $g_{ij}$ is expressed in the $u$ coordinates. Since the metrics are uniformly bounded away from degeneracy, the inverse matrices $(g_k)^{in}$ also converge strongly to $g^{in}$ in $C^{0,\alpha}$. Using the identity $\partial (g_k)^{in} = -(g_k)^{ia}(\partial (g_k)_{ab})(g_k)^{bn}$, the product of this convergence with the weak $L^q$ convergence of $\partial (g_k)_{ab}$ shows that $(g_k)^{in} \rightharpoonup g^{in}$ weakly in $W^{1,q}(B_R^+)$.

	To begin, let us consider the component $g_k^{nn}$. By Proposition \ref{thm:boundary_regularity}, it satisfies the local weak Neumann problem
	\begin{equation} \label{eq:weak_neum_gnn}
		\begin{cases}
			\Delta_{g_k} g_k^{nn} = F_k     & \text{ in } B_R^+,    \\
			\partial_{\nu_k} g_k^{nn} = S_k & \text{ on } \Sigma_R,
		\end{cases}
	\end{equation}
	with $F_k = B^{nn}(g_k, \partial g_k) - 2 \Ric(g_k)^{nn}$ and $S_k = -2(n-1)H_{g_k} g_k^{nn}$. We have $F_k$ uniformly bounded in $L^{q/2}(B_R^+)$ and $S_k$ uniformly
	bounded in $L^{\infty}(\Sigma_R)$. Thus, after passing to a subsequence, $F_k$ converges weakly in $L^{q/2}(B_R^+)$ to a limit $F_{\infty} \in L^{q/2}(B_R^+)$ and $S_k$ converges weak-* to a limit $S_{\infty} \in L^{\infty}(\Sigma_R)$.

	For any test function $\psi \in C_c^\infty(B_R^+ \cup \Sigma_R)$, the local weak Neumann problem for $g_k^{nn}$ is
	\begin{equation}
		\int_{B_R^+} \langle d g_k^{nn}, d\psi \rangle_{g_k} \sqrt{|g_k|} \, dx = - \int_{B_R^+} F_k \psi \sqrt{|g_k|} \, dx - \int_{\Sigma_R} S_k \psi \sqrt{|h_k|} \, dx'.
	\end{equation}
	Due to the strong convergence of $g_k \to g$ in $C^0(B_R^+)$ and $h_k \to h$ in $C^0(\Sigma_R)$, along with the weak convergence of $d g_k^{nn} \rightharpoonup d g^{nn}$ in $L^q(B_R^+)$, $F_k \rightharpoonup F_\infty$ in $L^{q/2}(B_R^+)$, and $S_k \stackrel{*}{\rightharpoonup} S_\infty$ in $L^\infty(\Sigma_R)$, we may pass to the limit as $k \to \infty$ to obtain
	\begin{equation}\label{eq:SF}
		\int_{B_R^+} \langle d g^{nn}, d\psi \rangle_{g} \sqrt{|g|} \, dx = - \int_{B_R^+} F_\infty \psi \sqrt{|g|} \, dx - \int_{\Sigma_R} S_\infty \psi \sqrt{|h|} \, dx'.
	\end{equation}
	This demonstrates that $F_\infty$ and $S_\infty$ solve the local weak Neumann problem for $g^{nn}$. Next, write $v_k = g_k^{nn} - g^{nn}$. We calculate
	\begin{align*}
		\int_{B_R^+} \langle d v_k, d\psi \rangle_{g_k} dV_{g_k} & = \int_{B_R^+} \langle d g_k^{nn}, d\psi \rangle_{g_k} dV_{g_k} - \int_{B_R^+} \langle d g^{nn}, d\psi \rangle_{g_k} dV_{g_k}               \\
		                                                         & = - \int_{B_R^+} F_k \psi \, dV_{g_k} - \int_{\Sigma_R} S_k \psi \, dS_{g_k} - \int_{B_R^+} \langle d g^{nn}, d\psi \rangle_{g_k} dV_{g_k}.
	\end{align*}
	Using \eqref{eq:SF}, we obtain
	\begin{align*}
		\int_{B_R^+} \langle d v_k, d\psi \rangle_{g_k} dV_{g_k} & = - \left( \int_{B_R^+} F_k \psi \, dV_{g_k} - \int_{B_R^+} F_\infty \psi \, dV_{g} \right)                                                  \\
		                                                         & \quad - \left( \int_{B_R^+} \langle d g^{nn}, d\psi \rangle_{g_k} dV_{g_k} - \int_{B_R^+} \langle d g^{nn}, d\psi \rangle_{g} dV_{g} \right) \\
		                                                         & \quad - \left( \int_{\Sigma_R} S_k \psi \, dS_{g_k} - \int_{\Sigma_R} S_\infty \psi \, dS_{g} \right),
	\end{align*}
	which implies that $v_k$ solves \eqref{eq:local_weak_neumann} with
	\begin{equation}\label{eq:F_def}
		\begin{aligned}
			(F, \psi) & = \int_{B_R^+} \left( F_k \sqrt{|g_k|} - F_\infty \sqrt{|g|} \right) \psi \, dx                                       \\
			          & \quad + \int_{B_R^+} \left( \sqrt{|g_k|} g_k^{ij} - \sqrt{|g|} g^{ij} \right) \partial_i g^{nn} \partial_j \psi \, dx
		\end{aligned}
	\end{equation}

	and

	\begin{equation}\label{eq:G_def}
		(G, \psi) = \int_{\Sigma_R} \left( S_k \sqrt{|h_k|} - S_\infty \sqrt{|h|} \right) \psi \, dx'.
	\end{equation}
	Using Corollary \ref{cor:local_neumann_w12}, we get
	\begin{equation}\label{eq:vk_quantitative_estimate}
		\|v_k\|_{H^1(B_r^+)} \le C \left( \|F\|_{H^{-1}(B_R^+)} + \|G\|_{H^{-1/2}(\Sigma_R)} + \|v_k\|_{L^2(B_R^+)} \right).
	\end{equation}
	Since $F_k \rightharpoonup F_{\infty}$ weakly in $L^{q/2}(B_R^+)$, after passing to a subsequence we have $F_k \to F_{\infty}$ strongly in $H^{-1}(B_R^+)$ due to the compact embedding $L^{q/2}(B_R^+) \hookrightarrow H^{-1}(B_R^+)$ (we choose $q$ large enough so that $q/2 > \frac{2n}{n+2}$ and apply the dual Sobolev embedding). Since $g_k \to g$ in $C^0(B_R^+)$, we conclude the first term on the right-hand side of \eqref{eq:F_def} tends to $0$ in $H^{-1}(B_R^+)$. Next, applying H\"older's inequality to the second term on the
	right-hand side of \eqref{eq:F_def} and noting that $\nabla g^{nn} \in L^q(B_R^+)$, we see this term tends to zero in $H^{-1}(B_R^+)$. Thus $\|F\|_{H^{-1}(B_R^+)} \to 0$ as $k \to \infty$. A parallel argument, using $S_k \stackrel{*}{\rightharpoonup} S_{\infty}$ in $L^\infty(\Sigma_R)$ and $h_k \to h$ in $C^0(\Sigma_R)$, shows that $\|G\|_{H^{-1/2}(\Sigma_R)} \to 0$ thanks to the compact embedding $L^p(\Sigma_R) \hookrightarrow H^{-1/2}(\Sigma_R)$ for $p \ge 2$. Finally, it's immediate from our hypothesis that $\|v_k\|_{L^2(B^+_R)} \to 0$.

	Since all terms on the right-hand side of \eqref{eq:vk_quantitative_estimate} vanish as $k \to \infty$, we have $v_k \to 0$ strongly in $H^1(B_r^+)$. Interpolation gives
	\begin{equation}
		\|v_k\|_{W^{1,q'}(B_{r}^+)} \le C \|v_k\|_{H^1(B_{r}^+)}^\theta \|v_k\|_{W^{1,q}(B_{r}^+)}^{1-\theta} \to 0.
	\end{equation}
	Thus, $g_k^{nn} \to g^{nn}$ strongly in $W^{1,q'}(B_r^+)$, $2 < q' < q$. Since $q$ was arbitrary, the convergence is strong in $W^{1,q}(B_r^+)$ for all $q > n$.

	Next, consider the components $g_k^{\ell n}$, $1 \le \ell \le n-1$. This case is parallel to the analysis of $g^{nn}$. Similarly to \eqref{eq:weak_neum_gnn}, we get
	\begin{equation} \label{eq:weak_neum_gelln}
		\begin{cases}
			\Delta_{g_k} g_k^{\ell n} = F_k^\ell     & \text{ in } B_R^+,    \\
			\partial_{\nu_k} g_k^{\ell n} = S_k^\ell & \text{ on } \Sigma_R,
		\end{cases}
	\end{equation}
	with $F_k^\ell = B^{\ell n}(g_k, \partial g_k) - 2 \Ric(g_k)^{\ell n}$ and $S_k^\ell = -(n-1)H_{g_k} g_k^{\ell n} + \frac{1}{2\sqrt{g_k^{nn}}} g_k^{\ell m} \partial_m g_k^{nn}$. In this case $F_k^{\ell}$ is uniformly bounded in $L^{q/2}(B_R^+)$ and $S_k^{\ell}$ is uniformly bounded in $L^{q}(\Sigma_R)$, any $q > n$. Therefore we have $F^{\ell}_{\infty}$ and $S_{\infty}^{\ell}$ so that $F_k^{\ell} \rightharpoonup F^{\ell}_{\infty}$ in $L^{q/2}(B_R^+)$ and $S^{\ell}_k \rightharpoonup S^{\ell}_{\infty}$ in $L^{q}(\Sigma_R)$.

	Write $w_k = g_k^{\ell n} - g^{\ell n}$. As before, $w_k$ satisfies a weak Neumann problem
	\begin{equation}\label{eq:weak_neum_wk}
		\begin{cases}
			\Delta_{g_k} w_k = F     & \text{ in } B_R^+,    \\
			\partial_{\nu_k} w_k = G & \text{ on } \Sigma_R,
		\end{cases}
	\end{equation}
	where
	\begin{equation}\label{eq:F_ell_def}
		\begin{aligned}
			(F, \psi) & = \int_{B_R^+} \left( F_k^\ell \sqrt{|g_k|} - F_\infty^\ell \sqrt{|g|} \right) \psi \, dx                                 \\
			          & \quad + \int_{B_R^+} \left( \sqrt{|g_k|} g_k^{ij} - \sqrt{|g|} g^{ij} \right) \partial_i g^{\ell n} \partial_j \psi \, dx
		\end{aligned}
	\end{equation}
	and
	\begin{equation}\label{eq:G_ell_def}
		(G, \psi) = \int_{\Sigma_R} \left( S_k^\ell \sqrt{|h_k|} - S_\infty^\ell \sqrt{|h|} \right) \psi \, dx'.
	\end{equation}

	Our hypotheses on $g_k$, $h_k$, and the weak convergence of $F^{\ell}_k$ and $S^{\ell}_k$ to their limits shows $\|F\|_{H^{-1}(B_R^+)} \to 0$ and $\|G\|_{H^{-1/2}(\Sigma_R)}\to 0$. Applying the local $H^1$ Neumann estimate \eqref{eq:local_neumann_w12_estimate} to $w_k$, we conclude $w_k \to 0$ strongly in $H^1(B_r^+)$. Interpolating with the uniform $W^{1,q}(B_R^+)$ bound ensures $g_k^{\ell n} \to g^{\ell n}$ strongly in $W^{1,q}(B_r^+)$.

	For the tangential components $(g_k)_{\alpha \beta}$, we can directly apply standard Dirichlet estimates. First, since $\{g_k\}$ are uniformly bounded in $W^{1,q}_{\mathrm{loc}}$ for any $q < \infty$, our hypotheses imply $\Delta_{h_k} (h_k)_{\alpha\beta} \in L^q(\Sigma_R)$, and elliptic regularity gives $(h_k)_{\alpha\beta} \in W^{2,q}_{\mathrm{loc}}(\Sigma_R)$, independent of $k$. Next, the uniform bound on $\Ric(g_k)$ and local elliptic estimates for the Dirichlet problem (Lemma \ref{thm:local_boundary_w2p_nonhom}) show that $(g_k)_{\alpha \beta}$ is uniformly bounded in $W^{2,q}_{\mathrm{loc}}(B_R^+)$. The compact embedding $W^{2,q} \hookrightarrow W^{1,q}$ implies these tangential components converge in $W^{1,q}_{\mathrm{loc}}(B_r^+)$.

	Finally, we establish convergence for the components $(g_k)_{in}$, $1 \leq i \leq n$. Dropping the index $k$ for the moment, note that $g^{ij} g_{jm} = \delta^i_m$ gives
	\begin{align*}
		g_{\alpha\beta} g^{\beta n} + g_{\alpha n} g^{nn} & = 0, \\
		g_{n\alpha} g^{\alpha n} + g_{nn} g^{nn}          & = 1,
	\end{align*}
	where $1 \le \alpha, \beta \le n-1$. Since $g^{nn}$ is strictly positive and bounded away from zero, we can solve for the remaining covariant components:
	\begin{align}
		g_{\alpha n} & = - \frac{1}{g^{nn}} g_{\alpha\beta} g^{\beta n}, \label{eq:alg_mixed}                                    \\
		g_{nn}       & = \frac{1}{g^{nn}} + \frac{1}{(g^{nn})^2} g_{\alpha\beta} g^{\alpha n} g^{\beta n}. \label{eq:alg_normal}
	\end{align}
	Because $W^{1,q}_{\mathrm{loc}}(B_r^+)$ forms a Banach algebra for $q > n$, it is closed under multiplication and division by functions bounded away from zero. These algebraic relations \eqref{eq:alg_mixed} and \eqref{eq:alg_normal} demonstrate that the strong convergence of $(g_k)_{\alpha\beta}$, $g_k^{\alpha n}$, and $g_k^{nn}$ fully determines the strong convergence of $(g_k)_{\alpha n}$ and $(g_k)_{nn}$. Therefore, $g_k \to g$ strongly in $W^{1,q}_{\mathrm{loc}}(B_r^+)$.
\end{proof}

\begin{proposition}\label{weaktostrong}
	Suppose $(M_k, g_k) \in \mathcal A^{bdd}$ is a sequence with $(M_k, g_k) \rightharpoonup (M, g)$ weakly in $W^{1,p}$ and $(\partial M_k, h_k) \rightharpoonup (\partial M, h)$ weakly in $W^{1,p}$, $p > n$. Then $(M_k, g_k) \to (M,g)$ subconverges strongly in $W^{1,q}$, any $q > n$.
\end{proposition}

\begin{proof}
	From the definition of weak convergence, the metrics $g_k$ are uniformly bounded in $W^{1,p}(M_k)$, $g \in W^{1,p}(M)$, $h_k$ are uniformly bounded in $W^{1,p}(\partial M_k)$, and $h \in W^{1,p}(\partial M)$. By Proposition \ref{thm:boundary_regularity}, the $g_k$ are uniformly bounded in $W^{1,q}(M_k)$, any $q > n$, and $g \in W^{1,q}(M)$. From \cite{MR1074481}, the $h_k$ subconverge to $h$ in the $C^{1,\alpha}$ topology, any $0 < \alpha < 1$. Moreover, the Sobolev embedding implies that $g_k \to g$ in the $C^{0,\alpha}$
	topology. We may also assume, after removing the first finitely many terms of $\{ g_k \}$, that each $g_k$ is pulled back to $M$ via a diffeomorphism and work on the fixed limit space $M$. We aim to construct harmonic atlases for $g_k$ that converge to a harmonic atlas for $g$ on fixed domains, which will allow us to use Proposition \ref{lem:local_strong_convergence}.

	Using the proof of Lemma \ref{thm:boundary_harmonic_existence}, cover a collar neighborhood of $\partial M$ with a finite collection of boundary harmonic coordinate charts $u_{\alpha}$, which are $W^{2,q}$ functions on their domains, any $q > n$. Write $v_{\alpha}$ for the corresponding finite atlas of $h$-harmonic charts on $\partial M$. Fix a point $z \in
		\partial M$ and choose such a chart (dropping the index), $u : U \to \Omega_R$ centered at $z \in U$, with corresponding boundary chart $v : U \cap \partial M \to \Sigma_R$, where $\Sigma_R = \overline{\Omega}_R \cap \{ x^n = 0 \}$. Here as in Lemma \ref{thm:boundary_harmonic_existence}, $\Omega_R$ is a smoothing of $B_R^+$, uniformly close to $B_R^+$ and with
	$\Sigma_R = B_R^+ \cap \{ x^n = 0\}$. By hypothesis and because of optimal regularity of harmonic coordinates, $(h_k)_{ij} \to h_{ij}$ in $C^{0, \alpha}(\Sigma_{R})$. In fact, as we noted above, we can assume that $h_k \to h$ in $C^{1,\alpha}$ on $\Sigma_{R}$. Next, with a slightly smaller radius $R_1$, solve $\Delta_{h_k} v_k^i = 0$ on $\Sigma_{R_1}$ with boundary values $v^i_k|_{\Gamma_{R_1}} = v^i$, where $\Gamma_{R_1} = \partial \Sigma_{R_1}$ (cf. Lemma \ref{thm:strong_solutions}). Then construct harmonic functions $u_k^i$ solving $\Delta_{g_k} u_k^i = 0$ on $\Omega_{R_1}$, and with
	boundary values $u_k^i = v_k^i$ on $\Sigma_{R_1}$ and $u_k^i = u^i$ on $\partial \Omega_{R_1} \setminus \Sigma_{R_1}$ (cf. Lemma \ref{thm:strong_solutions_c0}). By the maximum principle \cite[Theorem 9.1]{MR1814364}, the functions $u_k = (u_k^1, \ldots, u_k^n)$ are uniformly bounded in $C^{0}(\Omega_{R_1})$ independent of $k$, and we initially have $u_k \in W^{2,q}_{loc}(\Omega_{R_1})\cap C^0(\overline \Omega_{R_1})$. Returning to $v_k^i$, Lemma \ref{thm:strong_solutions} gives
	\begin{equation}\label{eq:boundary_control}
		\|v_k^i\|_{W^{2,q}(\Sigma_{R_1})} \le C \Big( \|v_k^i\|_{W^{2-1/q,q}(\Gamma_{R_1})} + \|v_k^i\|_{L^q(\Sigma_{R_1})} \Big).
	\end{equation}
	With $v_k^i = v^i$ on $\Gamma_{R_1}$ and the Alexandrov maximum principle, this shows $v_k^i$ is uniformly bounded in $W^{2,q}(\Sigma_{R_1})$. Thus, there is a subsequential limit $v_k^i \to \tilde v^i$ in $C^{1,\alpha}(\Sigma_{R_1})$ and weakly in $W^{2,q}(\Sigma_{R_1})$. Because $h_k \to h$ in $C^0(\Sigma_{R_1})$, $\tilde v^i$ satisfies the weak Dirichlet problem $\Delta_h \tilde v^i = 0$ with $\tilde v^i = v^i$ on $\Gamma_{R_1}$. 	By the weak maximum principle \cite[Theorem 8.1]{MR1814364}, we conclude $\tilde v^i = v^i$.

	Next, letting $V_k^i$ denote the piecewise boundary data defined by $V_k^i = v_k^i$ on $\Sigma_{R_1}$ and $V_k^i = u^i$ on $\partial \Omega_{R_1} \setminus \Sigma_{R_1}$, we estimate $u_k^i$ on a slightly smaller ball of radius $R' < R_1$.

	By Lemma \ref{thm:local_boundary_w2p_nonhom}, we have
	\begin{equation}
		\|u_k^i\|_{W^{2,q}(B_{R'}^+)} \le C \Big( \|v_k^i\|_{W^{2-1/q,q}(\Sigma_{R_1})} + \|u_k^i\|_{L^q(\Omega_{R_1})} \Big).
	\end{equation}
	Since $u^i_k$ is uniformly bounded in $C^0(\Omega_{R_1})$ and $v_k^i$ is uniformly bounded in $W^{2,q}(\Sigma_{R_1})$, this shows $u_k^i$ is uniformly bounded in $W^{2,q}(B_{R'}^+)$.
	Passing to a subsequence, we get that $u^i_k$ subconverges weakly in $W^{2,q}(B_{R'}^+)$. Taking note of the boundary values of $u^i_k$, using a diagonal argument and letting $R' \to R$, we obtain a limit $\tilde u^i \in W^{2,q}_{loc}(B_{R_1}^+) \cap C^0(\overline{B_{R_1}^+})$ with $u^i_k \rightharpoonup \tilde u^i$ in $W^{2,p}$ on compact subsets of $B^+_{R_1} \cup \Sigma_{R_1}$, and $\tilde u^i = u^i$ on $\partial B_{R_1}^+$. The weak convergence $g_k \rightharpoonup g$ implies $\tilde u^i$ is a weak solution of $\Delta_g \tilde u^i = 0$ on $B^+_{R_1}$, so the weak maximum principle implies $\tilde u^i = u^i$.

	Thus, for large enough $k$ we see that $u_k$ are boundary harmonic coordinates on $B^+_{R_2}$, any $R_2 < R_1$. This procedure can be done for every $u_{\alpha}$ in the atlas. In particular, the conditions of Proposition \ref{lem:local_strong_convergence} are satisfied on $B^+_{R_2}$, and we therefore conclude that $g_k \to g$ strongly in $W^{1,q}(B_{R_3})$ in each chart, for some smaller $R_3< R_2$, which implies $g_k \to g$ in the $W^{1,q}_{loc}$ topology on a collar neighborhood of $M$.

	The interior charts are easier to handle. The results here are standard, and stronger; one has weak $W^{2,q}$ and strong $C^{1,\alpha}$ metric convergence as a consequence of the arguments in \cite{MR1074481}.

	Because the covering of $M$ is finite, we may extract a single subsequence that converges strongly in $W^{1,q}$ on every chart in the finite subcover. We therefore obtain global strong $W^{1,q}$ convergence of $g_k$ to $g$ on $M$.
\end{proof}

\begin{lemma} \label{lem:anderson_w1p}
	Assume that a sequence of pointed Riemannian $n$-manifolds with boundary $(M_k, g_k, p_k)$, with $p_k \in \partial M_k$, converges to a limit $(M, g, p)$ in the pointed $W^{1,q}$ topology for all $q > n$. Write $(\partial M_k, h_k, p_k)$ for the boundary metrics and $(\partial M, h, p)$ for the boundary of the limit. Assume further that as $k \to \infty$:
	\begin{equation*}
		\|\Ric(g_k)\|_{L^\infty(M_k)} \to 0, \quad \|\Ric(h_k)\|_{L^\infty(\partial M_k)} \to 0, \quad \|H_k\|_{L^\infty(\partial M_k)} \to 0,
	\end{equation*}
	and
	\begin{equation*}
		i_{\partial M_k} \to \infty, \quad i_{b}(M_k) \to \infty.
	\end{equation*}
	Then $(M, g)$ is isometric to $\mathbb{R}^n_+$.
\end{lemma}

\begin{proof}
	This is essentially the same as the proof given in \cite[Lemma 3.2.2]{MR2096795}, but we must verify that the result holds good under our hypothesis of $W^{1,q}$ convergence, instead of the $C^{1,\alpha}$ convergence assumed there.

	First, as in Theorem \ref{precompact}, we have that $\Ric(g) = 0$, $\Ric_{\partial M}(h) = 0$, and $H = 0$, and thus $M$ is a smooth, Ricci-flat manifold with boundary. In fact, from \cite{MR1074481}, we know that $\partial M$ is isometric to $\R^{n-1}$. We can explicitly show that $(M,g)$ is
	isometric to $\mathbb{R}^n_+$ by analyzing the distance functions $\tau_k(x) = \dist_{g_k}(x, \partial M_k)$. Since $i_{b}(M_k) \to \infty$, the function $\tau_k$ is smooth on any fixed compact subset of $M_k$ for sufficiently large $k$. Let $H_k(t)$ denote the mean curvature of the level set $\Sigma_t = \{\tau_k = t\}$, and we note that $\Delta_k \tau_k = (n-1)H_k$. Just as in \cite[Lemma 3.2.2]{MR2096795}, we have that $H_k(t) \to 0$ uniformly on compact subsets of $[0, \infty)$.

	Because $g_k \to g$ strongly in $W^{1,q}_{loc}$ for any $q > n$, we have $\tau_k \to \tau$ uniformly, where $\tau$ is the limit distance function. Further, the leading coefficients of $\Delta_k$ are uniformly bounded in $C^\alpha$ and the first-order coefficients in $L^q$. To establish uniform
	local bounds, we fix a local boundary harmonic coordinate chart $u: U \to B_R^+$ (cf. Lemma \ref{thm:boundary_harmonic_existence}), where the boundary is the flat slice $\Sigma_R =\overline{B_R^+} \cap \{x^n = 0\}$. Local elliptic estimates (cf. Lemma \ref{thm:local_boundary_w2p_nonhom}) give a bound on a slightly smaller half-ball $B_r^+$,
	\begin{equation} \label{eq:tau_apriori}
		\|\tau_k\|_{W^{2,q}(B_r^+)} \leq C \left( \|\Delta_k \tau_k\|_{L^q(B_R^+)} + \|\tau_k\|_{W^{2-1/q,q}(\Sigma_R)} + \|\tau_k\|_{L^q(B_R^+)} \right).
	\end{equation}
	Since $\tau_k \equiv 0$ on $\Sigma_R$, $H_k \to 0$ in $L^\infty(B_R^+)$, and $\|\Delta_k \tau_k\|_{L^q(B_R^+)} = (n-1)\|H_k\|_{L^q(B_R^+)}$, this guarantees that $\tau_k$ is uniformly bounded in $W^{2,q}(B_r^+)$. Similar interior estimates hold away from the boundary. Consequently, $\tau_k$ subconverges weakly in $W^{2,q}_{loc}(B_R^+)$ and strongly in $C^{1,\alpha}_{loc}(B_R^+)$ to $\tau$.

	In local coordinates, the Laplacian is given by $\Delta_k \tau_k = g_k^{ij} \partial_{ij} \tau_k - g_k^{ij} (\Gamma_k)_{ij}^\ell \partial_\ell \tau_k$. Because $g_k \to g$ strongly in
	$W^{1,q}_{loc}(B_R^+)$, the leading coefficients $g_k^{ij}$ converge uniformly in $C^0(B_R^+)$ and the Christoffel symbols $(\Gamma_k)_{ij}^\ell$ converge strongly in $L^q_{loc}(B_R^+)
	$. Meanwhile, $\tau_k \rightharpoonup \tau$ weakly in $W^{2,q}_{loc}(B_R^+)$ and strongly in $C^{1,\alpha}_{loc}(B_R^+)$, meaning the second derivatives $\partial_{ij} \tau_k$ converge
	weakly in $L^q_{loc}(B_R^+)$ and the first derivatives $\partial_\ell \tau_k$ converge uniformly. Therefore, we may pass to the limit in each term of the coordinate expression to
	obtain $\Delta_k \tau_k \rightharpoonup \Delta_g \tau$ weakly in $L^q_{loc}(B_R^+)$. Thus, for every test
	function $\psi \in L^{q'}(B_r^+)$, $1/q + 1/q' = 1$, we get $\int_{B_r^+}\Delta_g \tau \psi dx = \lim_{k\to \infty}\int_{B_r^+} (n-1)H_k \psi dx = 0$, and therefore we conclude $\Delta_g \tau = 0$.

	Now the proof proceeds just as in \cite[Lemma 3.2.2]{MR2096795}. Because $\Delta_g \tau$ is $(n-1)$ times the mean curvature $H(\tau)$ of the level sets of $\tau$ in the limit manifold, we conclude $H(\tau) \equiv 0$ on any fixed interval $[0,t]$. Applying the Riccati equation reveals that the second fundamental form is identically zero, and $g$ splits isometrically as $[0, \infty) \times \partial M \cong \R^n_+$.
\end{proof}

\begin{remark}\label{rmk:stratified_convergence}
	The convergence properties of the boundary, extrinsic geometry of the boundary, and interior geometry of the manifold in Theorem \ref{precompact} can be refined as follows.  The intrinsic boundary metrics $h_k$ converge weakly in $W^{2,q}(\partial M)$ to $h_\infty$, for any $q > n$. The metrics $g_k$ converge weakly in $W^{2,q}_{\mathrm{loc}}(M \setminus \partial M)$ to $g_\infty$, for any $q > n$. The metrics $g_k$ converge strongly in $W^{1,q}(M)$ to $g_\infty$ globally up to the boundary, for any $q > n$. (This is just the statement of Theorem \ref{precompact}.) The first two statements already follow from \cite{MR1074481}; the point of our remark is that these convergence properties can be obtained from the same sequence. To show this, we create another refinement of the abstract convergence Theorem \ref{thm:abstract_convergence}, requiring harmonic radii at different regularity scales appropriate to the intrinsic boundary, interior, and collar neighborhood of the boundary (the definition of boundary harmonic radius in Definition \ref{def:boundary_harmonic_radius} is well-suited for this stratification). This changes the a priori regularity in Theorem \ref{precompact}, but the proof remains the same.
\end{remark}

\section{Stability theorems}\label{proof_of_T1}
We will use Theorem \ref{precompact} to prove two stability theorems regarding $3$-manifolds with boundary. As noted in the introduction, Theorem \ref{T1} below can be viewed, in a certain sense, as a generalization of the rigidity theorem of Cohn-Vossen (\cite{Cohn-Vossen, MR2261749}) and Pogorelov (\cite{MR0346714}). To see this, write $(\Sigma, h_{\Sigma})$ for a smooth, closed, oriented surface with Gauss curvature $K_{\Sigma} > 0$. By the solution of the Weyl problem (cf.\ \cite{MR2261749}), there exists a smooth isometric embedding $\Sigma \to \R^3$ whose image is unique up to rigid motion. Let $N \subset \R^3$ be the convex solid region bounded by this image, and let $g_{Euc}$ be the standard Euclidean metric on $N$. We establish conditions on a general 3-manifold $(M,g)$ that ensure it must be geometrically close to $(N, g_{Euc})$. We follow our approach that geometric hypotheses on $M$ should only involve uniform $L^{\infty}$ bounds on curvature quantities, although higher regularity versions can certainly be established and follow from standard techniques.
\begin{theorem}\label{T1}
	Suppose $N \subset \R^3$ is a smooth, closed, connected, compact, convex set. Write $\Sigma = \partial N$ and $h_{\Sigma}$ for the induced boundary metric on $\Sigma$.
	Let $(M_i, g_i)$ be a sequence of compact, oriented, simply connected Riemannian $3$-manifolds with connected boundaries. Write $h_i$ for the induced metrics on $\partial M_i$ and $K_i$ for their Gauss curvatures. Assume there exist constants $H_0 > 1$ and $K_0 > 1$ such that for all $i$:
	\begin{align*}
		1/H_0 \leq H_i \leq H_0 \quad \text{and} \quad 1/K_0 \leq K_i \leq K_0.
	\end{align*}
	If $\|\sec(g_i)\|_{L^\infty} \to 0$ uniformly and $(h_i, h_{\Sigma})_{GH} \to 0$ as $i \to \infty$, then for $i$ sufficiently large, there exist $C^{1,\alpha}$ diffeomorphisms $f_i : N \to M_i$ such that $f_i^*g_i \to g_{Euc}$ in the $C^\alpha$ topology for any $0 < \alpha < 1$.
\end{theorem}

\begin{proof}
	We are given a sequence $(M_i, g_i)$ satisfying the hypotheses of the theorem. We will first show that for large $i$, $(M_i,g_i) \in \mathcal{M}_+$.
	If $\lambda_1, \lambda_2$ are the principal curvatures of $\partial M_i$ at some point on $\partial M_i$ with tangent plane $\Pi$, then the Gauss equation gives, for large enough $i$,
	\begin{equation}\label{gauss_equation}
		\lambda_1\lambda_2 =  K_i - \sec_{M_i}(\Pi) > 0.
	\end{equation}
	Since $H_i = \frac{\lambda_1 + \lambda_2}{2} > 0$, we see that $\lambda_1$ and $\lambda_2$ are both positive. The upper bound on $H$ then implies that each $\lambda_k$ is uniformly bounded from above, while the inequality $K_i > 1/K_0$ shows that $\lambda_k$ is uniformly bounded from below. Thus, there is a constant $\widetilde H_0 > 0$ so that
	\begin{equation*}
		1/\widetilde H_0 < \lambda_k < \widetilde H_0.
	\end{equation*}
	In particular, each $M_i$ is uniformly convex when $i$ is large.

	Let us find an upper bound for $\diam(M_i)$. Myers' theorem implies that the intrinsic diameter of $\partial M_i$ is bounded from above. Let $k$ be a negative lower bound for $\sec(g_i)$. If $i$ is large enough, then $-k$ can be chosen small enough so that
	\begin{equation*}
		\frac{1}{\widetilde H_0\sqrt{-k}} > 1.
	\end{equation*}
	Thus, there exists a positive solution $t_0$ to the equation
	\begin{equation*}
		\coth(\sqrt{-k}t) = \frac{1}{\widetilde H_0\sqrt{-k}},
	\end{equation*}
	and by Lemma \ref{prop:focal_bound}, it follows that $\foc(\partial M_i) \leq t_0$. For any $p \in M_i$, there exists a length-minimizing geodesic from $p$ to $\partial M_i$ that meets $\partial M_i$ orthogonally. Since normal geodesics do not minimize distance to the boundary past the first focal point, it follows that $\dist(p, \partial M_i) \leq t_0$. Thus, $\diam(M_i)$ is bounded above by $\diam(\partial M_i) + 2t_0$. The Gauss-Bonnet theorem and the inequality $1/K_0 < K_i < K_0$ imply that $\vol_{n-1}(\partial M_i)$ is uniformly bounded below when $i$ is large enough.

	Therefore, eventually $(M_i, g_i) \in \mathcal{M}_+$ for appropriately defined constants. Theorem \ref{precompact} then shows that, after passing to a subsequence, $(M_i,g_i)$ converges in $W^{1,q}$, any $q > n$, to a limit $(M_{\infty},g_{\infty})$ with $g_{\infty} \in W^{1,q}(M_{\infty})$. Moreover, by Remark \ref{rmk:stratified_convergence}, the metrics $g_i$ converge in the weak $W^{2,q}_{\mathrm{loc}}$ or $C^{1,\alpha}$ topology on the interior $M_{\infty}\backslash \partial M_{\infty}$, and $h_i \to h_{\infty}$ in the weak $W^{2,q}$ (or $C^{1,\alpha}$) topology, any $0 < \alpha < 1$, where $h_i$ are the boundary metrics of $g_i$ and $h_{\infty}$ is the boundary metric of $g_{\infty}$. Since $h_i \to h_\Sigma$ in the Gromov-Hausdorff topology, we see that $(\partial M_{\infty},h_{\infty})$ and $(\Sigma,h_\Sigma)$ are isometric. Since $\partial M_{\infty}$ is orientable and admits a metric of positive curvature, it follows that $\partial M_{\infty}$ is diffeomorphic to $S^2$.

	In addition, we have $\sec(g_{\infty}) = 0$ (and thus $\Ric(g_{\infty}) = 0$) on $M_{\infty} \backslash \partial M_{\infty}$. Just as in Lemma \ref{lem:anderson_w1p}, the interior of $M_{\infty}$ is a smooth, flat Riemannian manifold. Since $M_{\infty}$ is simply connected, its interior $M_{\infty}\backslash \partial M_{\infty}$ inherits this simple connectivity (being homotopy equivalent to $M_{\infty}$). Therefore, we may construct the developing map, which is a smooth isometric immersion:
	\begin{align*}
		\xi: M_{\infty}\backslash \partial M_{\infty} \to \R^3.
	\end{align*}
	We now show that $\xi := (\xi^{\alpha})$ extends to the boundary with $W^{2,p}$ regularity. Writing $\delta$ for the Euclidean metric on $\R^3$, we have that $g_{\infty} = \xi^*\delta$, meaning $(D\xi)^T D\xi = g_{\infty}$. Evaluating the diagonal components of $g_{\infty}$ in any $W^{2,q}$ harmonic coordinate system $x$ gives $\sum_{\alpha=1}^3 (\partial_i \xi^\alpha)^2 = (g_{\infty})_{ii}$. Since $g_{\infty} \in C^{0,\alpha}(M_{\infty})$, we can certainly conclude that $\partial_i \xi^\alpha \in L^\infty(M_{\infty})$. Using again that $\xi$ is a local isometry, we turn to the Hessian equation Hess$(\xi^{\alpha}) = 0$, which gives
	\begin{align*}
		\frac{\partial^2 \xi^\alpha}{\partial x^i \partial x^j} = \Gamma_{ij}^k(g_{\infty}) \frac{\partial \xi^\alpha}{\partial x^k}.
	\end{align*}

	Since $\Gamma^k_{ij} \in L^q(M_{\infty})$  and $D\xi \in L^{\infty}(M_{\infty})$, we have that $\partial_i\partial_j \xi \in L^q(M_{\infty})$, any $q > n$. Consequently, $\xi$ belongs to $W^{2,q}(M_{\infty})$ up to the boundary. By applying the Sobolev embedding $W^{2,q} \hookrightarrow C^{1,\alpha}$ for $\alpha = 1 - 3/q$, we conclude that $\xi$ extends continuously to a $C^{1,\alpha}$ isometric immersion $M_{\infty} \to \R^3$, any $0 < \alpha < 1$.

	Let us show that the interior of $\xi(M_{\infty})$ is a convex region in $\mathbb{R}^3$. Fix $q_1, q_2 \in \xi( M_{\infty}\setminus \partial M_{\infty})$ and choose corresponding preimages $p_1, p_2$. Let $\tau_i(x) = \text{dist}_{g_i}(x, \partial M_i)$ denote the distance to the boundary in $M_i$. On a uniform collar neighborhood of $\partial M_i$, the Hessian of the distance function is given by $\text{Hess}(\tau_i)(X, X) = -g_i(S_i
		(X), X)$, where $S_i(X) = -\nabla_X \nu_i$ is the shape operator of the level sets with respect to the unit normal $\nu_i = \nabla\tau_i$. Because the boundaries are uniformly strictly convex for large $i$, $\text{Hess}(\tau_i)$ is negative semi-definite in this neighborhood.
	Let $\gamma_i : [0, 1] \to \Omega$ be the $g_i$-geodesic connecting $p_1$ and $p_2$, parametrized proportionally to $g_i$-arclength. Whenever $\gamma_i$ enters the collar neighborhood, the distance function evaluated along the geodesic, $f(s) = \tau_i(\gamma_i(s))$, satisfies $f''(s) = \text{Hess}(\tau_i)(\gamma_i', \gamma_i') \le 0$. This concavity implies that $\gamma_i$ cannot attain a local minimum distance to the boundary
	strictly inside the collar; therefore, $\tau_i(\gamma_i(s)) \ge \min(\tau_i(p_1), \tau_i(p_2))$. Since $p_1$ and $p_2$ are fixed strictly interior
	points, there exists a uniform constant $D > 0$ such that $\tau_i(\gamma_i(s)) \ge D$ for all large $i$. Thus, the geodesics $\gamma_i$ are uniformly confined to a compact region $\Omega \subset M_{\infty}\setminus \partial M_{\infty}$ upon which $g_i \to g_{\infty}$ in the $C^{1,\alpha}$ topology, for any $0 < \alpha < 1$.

	Let us show the $\gamma_i$ subconverge, with respect to the fixed metric $g_{\infty}$, in $C^{2,\alpha}([0,1])$, any $0 < \alpha <1$, to a limiting geodesic. So far, we have that $\{\gamma_i\}$ are uniformly bounded in $C^0([0,1])$. The choice of parametrization combined with uniform convergence $g_i \to g_{\infty}$ implies that $\dot{\gamma}$ is uniformly bounded in $C^{0}([0,1])$. Next, we consider the geodesic equation $\ddot{\gamma}_i^k = -(\Gamma_i)^k_{lm}(\gamma_i) \dot{\gamma}_i^l \dot{\gamma}_i^m$. Since the Christoffel symbols are uniformly bounded on $\Omega$, this shows that $\ddot{\gamma}_i$ is uniformly bounded in $C^0([0,1])$, making $\dot{\gamma}_i$ uniformly Lipschitz and thus uniformly bounded in $C^{0,\alpha}([0,1])$. Since the Christoffel symbols are also uniformly H\"older continuous, the right-hand side of the geodesic equation is therefore uniformly bounded in $C^{0,\alpha}([0,1])$. This implies that $\{\gamma_i\}$ is uniformly bounded in the $C^{2,\alpha}$ topology, and so a subsequence converges in the $C^{2,\alpha'}([0,1])$ topology to a limit curve $\gamma$. Because $0 < \alpha < 1$ was arbitrary, and passing to the subsequence without relabeling, we get $\gamma_i \to \gamma$ in $C^{2,\alpha}([0,1])$ for all $0 < \alpha < 1$. Therefore $\gamma$ is a $g_{\infty}$-geodesic, contained in $\Omega$, from $p_1$ to $p_2$. Therefore $\xi(\gamma)$ is strictly contained within the interior $\xi(M_{\infty})$. Because $g_{\infty}$ is a flat metric on a simply connected interior, its geodesics map to Euclidean straight lines under the local isometry $\xi$, meaning the straight line segment from $q_1$ to $q_2$ is contained in $\xi(M_{\infty} \backslash \partial M_{\infty})$, so that $\xi(M_{\infty})$ is convex.

	In particular, $\xi$ is globally injective, making it a $C^{1, \alpha}$ embedding. Composing the restriction of $\xi$ to $\partial M_{\infty}$ with a $C^{1,\alpha}$ isometry $(\Sigma,h_\Sigma) \to (\partial M_{\infty},h_{\infty})$, we obtain an isometric embedding of $(\Sigma,h_\Sigma)$ into $(\R^3,g_{Euc})$ whose image bounds a convex connected open set. A theorem of Pogorelov \cite[Thm 3.1.6]{MR0346714} implies that $\xi(\partial M_{\infty})$ differs from $\partial N$ by a rigid motion of $\R^3$. Composing $\xi$ by this rigid motion, we obtain the required $C^{1,\alpha}$ isometry between $M_{\infty}$ and $N$.
\end{proof}

We can quickly conclude a version of this using an $\epsilon$-$\delta$ framing.

\begin{corollary}\label{T1_eps_delta}
	Let $N \subset \mathbb{R}^3$ be a smooth, closed, connected, compact, convex set. Write $\Sigma = \partial N$ and $h_\Sigma$ for the induced boundary metric on $\Sigma$.
	Given constants $H_0 > 1$, $K_0 > 1$, for every $\epsilon > 0$ there exists $\delta = \delta(\epsilon, H_0, K_0) > 0$ such that the following holds.

	Suppose $(M,g)$ is a compact, oriented, simply connected Riemannian $3$-manifold with connected boundary. Write $h$ for the induced metric on $\partial M$, $K$ for the Gauss curvature of $h$, and $H$ for the mean curvature of $\partial M$. Assume that $H$ and $K$ satisfy the uniform bounds:
	\begin{align*}
		1/H_0 \leq H \leq H_0 \quad \text{and} \quad 1/K_0 \leq K \leq K_0.
	\end{align*}

	If the sectional curvature and boundary metric satisfy the conditions
	\begin{equation*}
		(h, h_{\Sigma})_{GH} < \delta, \qquad |\sec(g)| < \delta,
	\end{equation*}
	then there exists, for any $0 < \alpha < 1$, a $C^{1,\alpha}$ diffeomorphism $f:N \to M$ such that
	\begin{equation*}
		\|f^*g - g_{Euc}\|_{C^{\alpha}(N)} \leq \epsilon,
	\end{equation*}
	where $g_{Euc}$ is the standard Euclidean metric on $N$ and $(\cdot,\cdot)_{GH}$ is the Gromov-Hausdorff distance.
\end{corollary}

\begin{proof}
	This follows immediately from Theorem \ref{T1} by arguing by contradiction. If the corollary were false, there would exist an $\epsilon > 0$ and a sequence of manifolds $(M_i, g_i)$ satisfying the given bounds, with $|\sec(g_i)| \to 0$ and $(h_i, h_{\Sigma})_{GH} \to 0$ as $i \to \infty$, but for which no diffeomorphism $f_i: N \to M_i$ achieves $\|f_i^* g_i - g_{Euc}\|_{C^\alpha(N)} \leq \epsilon$. This directly contradicts the $C^\alpha$ convergence established in Theorem \ref{T1}.
\end{proof}

As a second application of Theorem \ref{precompact}, we present a geometric stability theorem motivated by Hopf's rigidity theorem (\cite{MR0040042}). Hopf's classic result states that any $C^3$ isometric immersion $i:S^2 \to \R^3$ of a topological sphere with constant mean curvature must have a round sphere as its image. Theorem \ref{conformal_stability} provides a quantitative, intrinsic stability version of this rigidity: a manifold that is approximately Euclidean and has a boundary with approximately constant mean curvature must be geometrically close to a Euclidean ball.

\begin{theorem}\label{conformal_stability}
	Given constants $S_0 > 0$, $D_0 > 0$, $C_0 > 0$, and $\alpha \in (0,1)$, for every $\epsilon > 0$ there exists $\delta = \delta(\epsilon, S_0, D_0, C_0, \alpha) > 0$ such that the following holds.

	Suppose $(M,g)$ is a compact, oriented, simply-connected Riemannian $3$-manifold with connected boundary satisfying $\chi(\partial M) = 2$, uniform bounds on the second fundamental form $|S| \leq S_0$, $\diam(M) \leq D_0$, and boundary metric $\|h\|_{C^{2,\alpha}(\partial M)} \leq C_0$. If there holds
	\begin{equation*}
		|\sec(g)| \leq \delta, \qquad \|H-1\|_{L^\infty} \leq \delta,
	\end{equation*}
	then there exists a $C^{1,\alpha}$ diffeomorphism $f: B \to M$ such that
	\begin{equation*}
		\|f^*g - g_{Euc}\|_{C^{\alpha}(B)} \leq \epsilon,
	\end{equation*}
	where $B \subset \R^3$ is the unit ball equipped with the standard Euclidean metric $g_{Euc}$.
\end{theorem}

\begin{proof}[Proof of Theorem \ref{conformal_stability}]
	We argue by contradiction. Suppose there exists an $\epsilon > 0$ and a sequence $(M_i, g_i)$ of compact, oriented Riemannian $3$-manifolds satisfying $\chi(\partial M_i) = 2$, $|S_i| \leq S_0$, $\diam(M_i) \leq D_0$, and $\|h_i\|_{C^{2,\alpha}(\partial M_i)} \leq C_0$, but with $|\sec(g_i)| \to 0$ and $\|H_i - 1\|_{L^\infty} \to 0$, such that no $C^{1,\alpha}$ diffeomorphism $f_i: B \to M_i$ achieves $\|f_i^*g_i - g_{Euc}\|_{C^\alpha} \leq \epsilon$.

	Let us first verify that $(M_i, g_i)$ eventually belongs to $\mathcal{M}_+$. Our hypotheses imply that the boundary Gauss curvature $K_i$ is bounded from above by some constant $K_{\max}$. By the Gauss-Bonnet theorem,
	\begin{equation*}
		4\pi = \int_{\partial M_i} K_i \, d\vol \leq K_{\max} \vol_{n-1}(\partial M_i),
	\end{equation*}
	which guarantees a uniform lower bound on $\vol_{n-1}(\partial M_i)$. Since $H_i \to 1$, the mean curvature is strictly bounded away from zero for large $i$. Hence, $(M_i, g_i) \in \mathcal{M}_+$. By Theorem \ref{precompact}, $(M_i,g_i)$ subconverges in the $W^{1,q}$ topology, any $q > n$, to a limit $(M_{\infty}, g_{\infty})$ satisfying
	\begin{equation}\label{eq:hopf_limit}
		H  = 1, \qquad \sec(g_{\infty}) = 0.
	\end{equation}
	From the uniform boundedness of $h_i$ in $C^{2,\alpha}(\partial M_i)$, we may assume the intrinsic boundary metric of the limit $h_{\infty} \in C^{2,\alpha}$ and that $h_i \to h_{\infty}$ in $C^{2,\alpha'}$, any $0 < \alpha' < \alpha$. Combined with \eqref{eq:hopf_limit}, elliptic regularity implies that $g_{\infty} \in C^{2,\alpha}(M_{\infty})$, any $0 < \alpha < 1$.

	As in \ref{T1}, we consider the developing map $\xi: M_{\infty} \to \R^3$. This time, $g \in C^{2,\alpha}(M_{\infty})$ implies $\xi \in C^{3, \alpha}(M_{\infty})$. Restricting $\xi$ to $\partial M_{\infty}$ gives a $C^3$ isometric immersion of a surface of genus $0$ with constant mean curvature $H = 1$. Thus,  Hopf's rigidity theorem applies and the image $\xi(M_{\infty}) = B.$ This implies that $(M_i, g_i)$ converges to $(B, g_{Euc})$ in the strong $C^\alpha$ topology, which contradicts our initial assumption and completes the proof.
\end{proof}

\bibliographystyle{amsplain}
\bibliography{bibliography}
\end{document}